\documentclass[10pt]{amsart}
\usepackage{verbatim}
\usepackage{eucal,url,amssymb,stmaryrd,enumerate,amscd,}
\usepackage[pagebackref,colorlinks=true]{hyperref}  
\usepackage{amsfonts}
\usepackage{amsmath,amsthm,amssymb,amscd,enumerate,eucal,url,stmaryrd}

\numberwithin{equation}{section}
\newtheorem{thrm}{Theorem}[section]
\newtheorem{lemma}[thrm]{Lemma}
\newtheorem{prop}[thrm]{Proposition}

\newtheorem{dfn}[thrm]{Definition}

\newtheorem{rmrk}[thrm]{Remark}

\newtheorem{conv}[thrm]{Convention}
\begin{document}

\begin{abstract}
{\ A complete solution to the quaternionic contact Yamabe problem on the
seven dimensional sphere is given. Extremals for the Sobolev inequality on
the seven dimensional Heisenberg group are explicitly described and the
best constant in the $L^2$ Folland-Stein embedding theorem is determined.}
\end{abstract}

\title[Extremals for the Sobolev inequality]{Extremals for the Sobolev
inequality on the seven dimensional quaternionic Heisenberg group and the
quaternionic contact Yamabe problem}
\date{\today }
\thanks{This project has been funded in part by the National Academy of
Sciences under the [Collaboration in Basic Science and Engineering Program 1
Twinning Program] supported by Contract No. INT-0002341 from the National
Science Foundation. The contents of this publication do not necessarily
reflect the views or policies of the National Academy of Sciences or the
National Science Foundation, nor does mention of trade names, commercial
products or organizations imply endorsement by the National Academy of
Sciences or the National Science Foundation.}
\author{Stefan Ivanov}
\address[Stefan Ivanov]{University of Sofia, Faculty of Mathematics and
Informatics, blvd. James Bourchier 5, 1164, Sofia, Bulgaria}
\email{ivanovsp@fmi.uni-sofia.bg}
\author{Ivan Minchev}
\address[Ivan Minchev]{University of Sofia\\
Sofia, Bulgaria\\
and Institut f\"ur Mathematik, Humboldt Universit\"at zu Berlin\\
Unter den Linden~6, Berlin~D-10099, Germany}
\email{minchevim@yahoo.com}
\email{minchev@fmi.uni-sofia.bg}
\author{Dimiter Vassilev}
\address[Dimiter Vassilev]{ Department of Mathematics and Statistics\\
University of New Mexico\\
Albuquerque, New Mexico, 87131-0001\\
and\\
University of California, Riverside\\
Riverside, CA 92521}
\email{vassilev@math.unm.edu}
\maketitle
\tableofcontents



\setcounter{tocdepth}{2}

\section{Introduction}

It is well known that the sphere at infinity of a any non-compact symmetric
space $M$ of rank one carries a natural Carnot-Carath\'eodory structure, see
\cite{M,P}. A quaternionic contact (qc) structure, \cite{Biq1,Biq2}, appears
naturally as the conformal boundary at infinity of the quaternionic
hyperbolic space. {In this paper, following Biquard, a quaternionic contact
structure (\emph{qc structure}) on a real (4n+3)-dimensional manifold $M$ is
a codimension three distribution $H$ locally given as the kernel of a $%
\mathbb{R}^3$-valued one-form $\eta=(\eta_1,\eta_2,\eta_3)$, such that, the
three two-forms $d\eta_i|_H$ are the fundamental forms of a quaternionic
structure on $H$. This means that there exists a Riemannian metric $g$ on $H$
and three local almost complex structures $I_i$ on $H$ satisfying the
commutation relations of the imaginary quaternions, $I_1I_2I_3=-1$, such
that, $d\eta_i|_H=2g(I_i.,.)$ . The 1-form $\eta$ is determined up to a
conformal factor and the action of $SO(3)$ on $\mathbb{R}^3$, and therefore $%
H$ is equipped with a conformal class $[g]$ of Riemannian metrics and a
2-sphere bundle of almost complex structures, the quaternionic bundle $%
\mathbb{Q}$. The 2-sphere bundle of one forms determines uniquely the
associated metric and a conformal change of the metric is equivalent to a
conformal change of the one forms. To every metric in the fixed conformal
class one can associate a linear connection preserving the qc structure, see
\cite{Biq1}, which we shall call the Biquard connection. }

If the first Pontrijagin class of $M$ vanishes then the 2-sphere bundle of $%
\mathbb{R}^3$-valued 1-forms is trivial \cite{AK}, i.e. there is a globally
defined 3-contact form $\eta$ that anihilates $H$, we denote the
corresponding QC manifold $(M,\eta)$. In this case the 2-sphere of
associated almost complex structures is also globally defined on $H$.

Examples of QC manifolds are given in \cite{Biq1,Biq2,IMV,D1}. In
particular, any totally umbilic hypersurface of a quaternionic K\"ahler or
hyperK\"ahler manifold carries such a structure \cite{IMV}. A basic example
is provided by any 3-Sasakian manifold which can be defined as a $(4n+3)$%
-dimensional Riemannian manifold whose Riemannian cone is a hyperK\"ahler
manifold. It was shown in \cite{IMV} that the torsion endomorphism of the
Biquard connection is the obstruction for a given qc-structure to be locally
3-Sasakian, up to a multiplication with a constant factor and a $SO(3)$%
-matrix.

{\ For a fixed metric in the conformal class of metrics on the horizontal
space one associates the scalar curvature of the associated Biquard
connection, called the qc-scalar curvature.} Guided by the real (Riemannian)
and complex (CR) cases, the quaternionic contact Yamabe problem is: \emph{%
given a compact QC manifold $(M,\eta)$, find a conformal 3-contact form for
which the qc-scalar curvature is constant.}

In the present paper
we provide a solution of this problem on the seven dimensional sphere
equipped with its natural quaternionic contact structure. The question
reduces to the solvability of the Yamabe equation
\eqref{e:conf change
scalar curv}. Taking the conformal factor in the form $\bar\eta=u^{4/(Q-2)}%
\eta$, $Q=4n+6$, turns \eqref{e:conf
change scalar curv} into the equation
\begin{equation*}
\mathcal{L} u\ \equiv\ 4\frac {Q+2}{Q-2}\ \triangle u -\ u\, Scal \ =\ - \
u^{2^*-1}\ \overline{Scal},
\end{equation*}
where $\triangle $ is the horizontal sub-Laplacian, $\triangle h\ =\
tr^g(\nabla dh)$, $Scal$ and $\overline{Scal}$ are the qc-scalar curvatures
correspondingly of $(M,\, \eta)$ and $(M, \, \bar\eta)$, and $2^* = \frac {2Q%
}{Q-2},$ with $Q=4n+6$--the homogeneous dimension. {On a compact
quaternionic contact manifold $M$ with a fixed conformal class $[\eta]$ the
Yamabe equation characterizes the non-negative extremals of the Yamabe
functional defined by
\begin{equation*}
\Upsilon (u)\ =\ \int_M\Bigl(4\frac {Q+2}{Q-2}\ \lvert \nabla u \rvert^2\ +\
\text{Scal}\, u^2\Bigr) dv_g,\qquad \int_M u^{2^*}\, dv_g \ =\ 1, \ u>0.
\end{equation*}
{\ Considering $M$ equipped with a fixed qc structure, hence, a conformal
class $[\eta]$,} the Yamabe constant is defined as the infimum
\begin{equation*}
\lambda(M)\ \equiv \ \lambda(M, [\eta])\ =\ \inf \{ \Upsilon (u) :\ \int_M
u^{2^*}\, dv_g \ =\ 1, \ u>0\}.
\end{equation*}
{\ Here $dv_g$ denotes the Riemannian volume form of the Riemannian metric
on $M$ extending in a natural way the horizontal metric associated to $\eta$.%
} }

When the Yamabe constant $\lambda(M)$ is less than that of the quaternionic
sphere with its standard qc structure the existence of solutions to the
quaternionic contact Yamabe problem is shown in \cite{Wei}, see also \cite%
{JL2}. We consider the Yamabe problem on the {\ standard unit $(4n+3)$%
-dimensional quaternionic sphere. The standard 3-Sasaki structure on the
sphere is a qc-Einstein structure $\tilde\eta$ having constant qc-scalar
curvature $\widetilde{\text{Scal}}=16n(n+2)$. Its images under conformal
quaternionic contact automorphism have again constant qc-scalar curvature.}
In \cite{IMV} we conjectured that these are the only solutions to the Yamabe
problem on the quaternionic sphere. The purpose of this paper is to prove
this conjecture when the dimension is equal to seven, i.e., $n=1$.

\begin{thrm}
\label{t:Yamabe} Let $\tilde\eta=\frac{1}{2h}\eta$ be a conformal
deformation of the standard qc-structure $\tilde\eta$ on the quaternionic
unit sphere $S^{7}$. If $\eta$ has constant qc-scalar curvature, then up to
a multiplicative constant $\eta$ is obtained from $\tilde\eta$ by a
conformal quaternionic contact automorphism. In particular, $\lambda(S^7)=
48\, (4\pi)^{1/5}$ and this minimum value is achieved only by $\tilde\eta$
and its images under conformal quaternionic contact automorphisms.
\end{thrm}

In \cite{IMV} a weaker result was shown, namely the conclusion holds (in all
dimensions) provided the vertical space of $\eta$ is integrable. {\ We
recall the definition of conformal quaternionic contact automorphism in
Definition \ref{d:3-ctct auto}.}

Another motivation for studying the Yamabe equation comes from its
connection with the determination of the norm and extremals in a relevant
Sobolev-type embedding on the quaternionic Heisenberg group $\boldsymbol{G\,(%
\mathbb{H})}$, \cite{GV} and \cite{Va1} and \cite{Va2}. As it is well known,
the Yamabe equation is essentially the Euler-Lagrange equation of the
extremals for the $L^2$ case of such embedding results. In the considered
setting we have the following Theorem due to Folland and Stein \cite{FS}.

\begin{thrm}[Folland and Stein]
\label{T:Folland and Stein} Let $\Omega \subset \boldsymbol{G}$ be an open
set in a Carnot group $\boldsymbol{G}$ of homogeneous dimension $Q$ {\ and
Haar measure $dH$}. For any $1<p<Q$ there exists $S_p=S_p(\boldsymbol{G})>0$
such that for $u\in C^\infty_o(\Omega)$
\begin{equation}  \label{FS}
\left(\int_\Omega\ |u|^{p^*}\ dH(g)\right)^{1/p^*} \leq\ S_p\
\left(\int_\Omega |Xu|^p\ dH(g)\right)^{1/p},
\end{equation}
where $|Xu|=\sum_{j=1}^m |X_ju|^2$ with $X_1,\dots, X_m$ denoting a basis of
the first layer of $\boldsymbol{G}$ {and $p^*= \frac {pQ}{Q-p}$.}
\end{thrm}

\noindent Let $S_p$ be the best constant in the Folland-Stein inequality,
i.e., the smallest constant for which \eqref{FS} holds. The {\ second}
result of this paper is the following Theorem, which determines the
extremals and the best constant in Theorem \ref{T:Folland and Stein} when $%
p=2$ for the case of the seven dimensional quaternionic Heisenberg group $%
\boldsymbol{G\,(\mathbb{H})}$. {\ As a manifold $\boldsymbol{G\,(\mathbb{H})}
\ =\mathbb{H}\times\text {Im}\, \mathbb{H}$ with the group law given by
\begin{equation*}
( q^{\prime }, \omega^{\prime })\ =\ (q_o, \omega_o)\circ(q, \omega)\ =\
(q_o\ +\ q, \omega\ +\ \omega_o\ + \ 2\ \text {Im}\ q_o\, \bar q),
\end{equation*}
\noindent where $q,\ q_o\in\mathbb{H}$ and $\omega, \omega_o\in \text {Im}\,
\mathbb{H}$. The standard quaternionic contact(qc) structure is defined by
the left-invariant quaternionic contact form $\tilde\Theta\ =\
(\tilde\Theta_1,\ \tilde\Theta_2, \ \tilde\Theta_3)\ =\ \frac 12\ (d\omega \
- \ q^{\prime }\cdot d\bar q^{\prime }\ + \ dq^{\prime }\, \cdot\bar
q^{\prime })$, where $.$ denotes the quaternion multiplication.}\ 

\begin{thrm}
\label{t:FS} Let $\boldsymbol{G\,(\mathbb{H})} \ =\ \mathbb{H}\times\text {Im%
}\, \mathbb{H}$ be the seven dimensional quaternionic Heisenberg group. The
best constant in the $L^2$ Folland-Stein embedding theorem is
\begin{equation*}  \label{best}
S_2\ =\ \frac{2\sqrt{3}}{\pi^{3/5}}
\end{equation*}
\noindent An extremal is given by the function

\begin{equation*}
v \ =\ \frac{2^{11}\sqrt{3}}{\pi^{3/5}}[(1+\lvert q \rvert^2)^2\ +\ \lvert
\omega \rvert^2]^{-2}, \ (q,\omega)\in \boldsymbol{G\,(\mathbb{H})}
\end{equation*}

\noindent Any other non-negative extremal is obtained from $v$ by
translations \eqref{e:translation} and dilations \eqref{e:scaling}.
\end{thrm}

{\ Our result confirms the Conjecture made after \cite[Theorem 1.1]{GV}. In
\cite[Theorem 1.6]{GV} the above Theorem is proved in all dimensions, but
with the assumption of partial-symmetry. Here with a completely different
method from \cite{GV} we show that the symmetry assumption is superfluous in
the case of the first quaternionic Heisenberg group. On the other hand, in
\cite{IMV} we proved Theorem~\ref{t:Yamabe} in all dimensions, but with the
'extra' assumption of the integrability of the vertical distribution.
In the present paper we remove the 'extra' integrability assumption in
dimension seven. A key step in present result is the establishment of a
suitable divergence formula, Theorem \ref{t:div formula}, see \cite{JL3} for
the CR case and \cite{Ob}, \cite{LP} for the Riemannian case. With the help
of this divergence formula we show that the 'new' structure is also
qc-Einstein, thus we reduce the Yamabe problem on the 7-sphere from solving
the non-linear Yamabe equation to a geometrical system of differential
equations describing the qc-Einstein structures conformal to the standard
one. Invoking the (quaternionic) Cayley transform, which is a contact
conformal diffeomorphism, \cite{IMV}, we turn the question to the
corresponding system on the quaternionic Heisenberg group. On the latter all
global solutions are explicitly described in \cite{IMV} and this allows us
to conclude the proof of our results.}

\begin{rmrk}
\label{rem111} {\ With the left invariant basis of Theorem~\ref{t:FS} the
Heisenberg group $\boldsymbol{G\,(\mathbb{H})}$ is not a group of Heisenberg
type. If we consider $\boldsymbol{G\,(\mathbb{H})}$ as a group of Heisenberg
type then the best constant in the $L^2$ Folland-Stein embedding theorem is,
cf. \cite[Theorem 1.6]{GV},
\begin{equation*}
S_2\ =\ \frac{15^{1/10}}{\pi^{2/5}\,2\,\sqrt{2}}.
\end{equation*}
\noindent and an extremal is given by the function
\begin{equation*}
F(q,\omega)\ =\ \gamma\ \left[(1 + |q|^2)^2\ +\ 16 |\omega|^2)\right]^{-2},
\ (q,\omega)\in \boldsymbol{G\,(\mathbb{H})}
\end{equation*}
where
\begin{equation*}
\gamma\ =\ 32\,\pi^{-17/50}\,2^{1/5}\,15^{2/5}.
\end{equation*}%
}
\end{rmrk}

\textbf{Organization of the paper.} The paper uses some results from \cite%
{IMV}. In order to make the present paper self-contained, in Section \ref%
{s:review} we give a review of the notion of a quaternionic contact
structure and collect formulas and results from \cite{IMV} that will be used
in the subsequent sections.

Section \ref{s:conf transf} and \ref{s:div formulas} are of technical
nature. In the former we find some transformations formulas for relevant
tensors, while in the latter we prove certain divergence formulas. The key
result is Theorem \ref{t:div formula}, with the help of which in the last
Section we prove the main Theorems.


\begin{conv}
\label{conven} We use the following conventions:
\begin{itemize}
\item
 $\{e_1,\dots,e_{4n}\}$ denotes an orthonormal basis of the
horizontal space $H$.
\item The summation convention over repeated vectors from the basis $%
\{e_1,\dots,e_{4n}\}$ will be used. For example, for a
(0,4)-tensor $P$, the formula $k=P(e_b,e_a,e_a,e_b)$ means
\begin{equation*}
k=\sum_{a,b=1}^{4n}P(e_b,e_a,e_a,e_b).
\end{equation*}
\item The triple $(i,j,k)$ denotes any cyclic permutation of
$(1,2,3)$.
\end{itemize}
\end{conv}

\textbf{Acknowledgements} S.Ivanov is visiting Max-Plank-Institut f\"ur
Mathematics, Bonn. S.I. thanks MPIM, Bonn for providing the support and an
excellent research environment during the final stages of the paper. S.I. is
a Senior Associate to the Abdus Salam ICTP. I.Minchev is a member of the
Junior Research Group "Special Geometries in Mathematical Physics" founded
by the Volkswagen Foundation. The authors would like to thank The National
Academies for the financial support and University of California, Riverside
and University of Sofia for hosting the respective visits of the authors.

{\ The authors would like to thank the referee for remarks making the
exposition clearer and spotting several typos in the paper. }

\section{Quaternionic contact manifolds}

\label{s:review} In this section we will briefly review the basic notions of
quaternionic contact geometry and recall some results from \cite{Biq1} and
\cite{IMV}.

For the purposes of this paper, a quaternionic contact (QC)
manifold $(M, g, \mathbb{Q})$ is a $4n+3$ dimensional manifold $M$
with a codimension three distribution $H$ equipped with a metric
$g$ and an Sp(n)Sp(1) structure, i.e., we have
\begin{enumerate}
\item[i)] a 2-sphere bundle $\mathbb{Q}$ over $M$ of almost
complex structures, such that, we have $\mathbb{Q}=
\{aI_1+bI_2+cI_3:\ a^2+b^2+c^2=1 \}$, where the almost complex
structures $I_s\,:H \rightarrow H,\quad I_s^2\ =\ -1, \quad s\ =\
1,\ 2,\ 3,$ satisfy the commutation relations of the imaginary
quaternions $I_1I_2=-I_2I_1=I_3$; \item[ii)] $H$ is the kernel of
a 1-form $\eta=(\eta_1,\eta_2,\eta_3)$ with values in
$\mathbb{R}^3$ and the following compatibility condition holds
\begin{equation*}  \label{con1}
2g(I_sX,Y)\ =\ d\eta_s(X,Y), \quad s=1,2,3, \quad X,Y\in H.
\end{equation*}
\end{enumerate}

Correspondingly, given a quaternionic contact manifold we shall denote with $%
\eta$ any associated contact form. The associated contact form is determined
up to an $SO(3)$-action, namely if $\Psi\in SO(3)$ with smooth functions as
entries then $\Psi\eta$ is again a contact form satisfying the above
compatibility condition (rotating also the almost complex structures) . On
the other hand, if we consider the conformal class $[g]$, the associated
contact forms are determined up to a multiplication with a positive function
$\mu$ and an $SO(3)$-action, namely if $\Psi\in SO(3)$ then $\mu\Psi\eta$ is
a contact form associated with a metric in the conformal class $[g]$.

We shall denote with $(M, \eta)$ a QC manifold with a fixed globally defined
contact form. A special phenomena here, noted in \cite{Biq1}, is that the
3-contact form $\eta$ determines the quaternionic structure and the metric
on the horizontal bundle in a unique way.

A QC manifold $(M, \bar g,\mathbb{Q} )$ is called conformal to $(M, g,%
\mathbb{Q} )$ if $\bar g\in [g]$. In that case, if $\bar\eta$ is a
corresponding associated one-form with complex structures $\bar I_s$, $%
s=1,2,3,$ we have $\bar\eta\ =\ \mu\, \Psi\,\eta$ for some $\Psi\in SO(3)$
with smooth functions as entries and a positive function $\mu$. In
particular, starting with a QC manifold $(M, \eta)$ and defining $\bar\eta\
=\ \mu\, \eta$ we obtain a QC manifold $(M, \bar\eta)$ conformal to the
original one.

\begin{dfn}
\label{d:3-ctct auto} A diffeomorphism $\phi$ of a QC manifold $(M,[g],%
\mathbb{Q})$ is called a \emph{conformal quaternionic contact automorphism
(conformal qc-automorphism)} if $\phi$ preserves the QC structure, i.e.
\begin{equation*}
\phi^*\eta=\mu\Psi\cdot\eta,
\end{equation*}
for some positive smooth function $\mu$ and some matrix $\Psi\in SO(3)$ with
smooth functions as entries and $\eta=(\eta_1,\eta_2,\eta_3)^t$ is a local
1-form considered as a column vector of three one forms as entries.
\end{dfn}

Any endomorphism $\Psi$ of $H$ can be decomposed with respect to the
quaternionic structure $(\mathbb{Q},g)$ uniquely into $Sp(n)$-invariant
parts as follows \hspace{2mm} 
$\Psi=\Psi^{+++}+\Psi^{+--}+\Psi^{-+-}+\Psi^{--+}, $ 
\hspace{2mm} where $\Psi^{+++}$ commutes with all three $I_i$, $\Psi^{+--}$
commutes with $I_1$ and anti-commutes with the other two and etc. \noindent
The two $Sp(n)Sp(1)$-invariant components are given by
\begin{equation}
{\label{New21}} \Psi_{[3]}=\Psi^{+++}, \qquad
\Psi_{[-1]}=\Psi^{+--}+\Psi^{-+-}+\Psi^{--+}.
\end{equation}
\noindent Denoting the corresponding (0,2) tensor via $g$ by the same letter
one sees that the $Sp(n)Sp(1)$-invariant components are the projections on
the eigenspaces of the Casimir operator
\begin{equation}  \label{e:cross}
\dag \ =\ I_1\otimes I_1\ +\ I_2\otimes I_2\ +\ I_3\otimes I_3
\end{equation}
corresponding, respectively, to the eigenvalues $3$ and $-1$, see \cite{CSal}%
. If $n=1$ then the space of symmetric endomorphisms commuting with all $%
I_i, i=1,2,3$ is 1-dimensional, i.e. the [3]-component of any symmetric
endomorphism $\Psi$ on $H$ is proportional to the identity, $\Psi_{[3]}=%
\frac{tr\, (\Psi)}{4}Id_{|H}$.

On a quaternionic contact manifold there exists a canonical connection
defined in \cite{Biq1} when the dimension $(4n+3)>7$, and in \cite{D} in the
7-dimensional case.

\begin{thrm}
\cite{Biq1}\label{biqcon} {Let $(M, g,\mathbb{Q})$ be a quaternionic contact
manifold} of dimension $4n+3>7$ and a fixed metric $g$ on $H$ in the
conformal class $[g]$. Then there exists a unique connection $\nabla$ with
torsion $T$ on $M^{4n+3}$ and a unique supplementary subspace $V$ to $H$ in $%
TM$, such that:

\begin{enumerate}[i)]
\item $\nabla$ preserves the decomposition $H\oplus V$ and
the metric $g$;
\item for $X,Y\in H$, one has $T(X,Y)=-[X,Y]_{|V}$;
\item
$\nabla$ preserves the $Sp(n)Sp(1)$-structure on $H$, i.e., $\nabla g\ = \
0$ and $\nabla \mathbb{Q}\subset\mathbb{Q}$;
\item for $\xi\in V$, the endomorphism $T(\xi,.)_{|H}$
of $H$ lies in $(sp(n)\oplus sp(1))^{\bot}\subset gl(4n)$;
\item the connection on $V$ is induced by the natural identification $\varphi$ of
$V$ with the subspace $sp(1)$ of the endomorphisms of $H$, i.e.
$\nabla\varphi=0$.
\end{enumerate}
\end{thrm}

We shall call the above connection \emph{the Biquard connection}.
Biquard \cite{Biq1} also described the supplementary subspace $V$
explicitly, namely, {locally }$V$ is generated by vector fields $%
\{\xi_1,\xi_2,\xi_3\}$, such that
\begin{equation}  \label{bi1}
\begin{aligned} \eta_s(\xi_k)=\delta_{sk}, \qquad (\xi_s\lrcorner
d\eta_s)_{|H}=0,\\ (\xi_s\lrcorner d\eta_k)_{|H}=-(\xi_k\lrcorner
d\eta_s)_{|H}. \end{aligned}
\end{equation}
The vector fields $\xi_1,\xi_2,\xi_3$ are called Reeb vector fields or
fundamental vector fields.

{\ If the dimension of $M$ is seven, the conditions \eqref{bi1} do not
always hold. Duchemin shows in \cite{D} that if we assume, in addition, the
existence of Reeb vector fields as in \eqref{bi1}, then Theorem~\ref{biqcon}
holds. Henceforth, by a qc structure in dimension $7$ we shall 
mean a qc structure satisfying \eqref{bi1}.}

Notice that equations \eqref{bi1} are invariant under the natural $SO(3)$
action. Using the triple of Reeb vector fields we extend $g$ to a metric on $%
M$ by requiring 
$span\{\xi_1,\xi_2,\xi_3\}=V\perp H \text{ and } g(\xi_s,\xi_k)=\delta_{sk}.
$ 
\hspace{2mm} \noindent The extended metric does not depend on the action of $%
SO(3)$ on $V$, but it changes in an obvious manner if $\eta$ is multiplied
by a conformal factor. Clearly, the Biquard connection preserves the
extended metric on $TM, \nabla g=0$. We shall also extend the quaternionic
structure by setting $I_{s|V}=0$. The fundamental 2-forms $\omega_i, i=1,2,3$
of the quaternionic structure $Q$ are defined by
\begin{equation}  \label{thirteen}
2\omega_{i|H}\ =\ \, d\eta_{i|H},\qquad \xi\lrcorner\omega_i=0,\quad \xi\in
V.
\end{equation}
Due to \eqref{thirteen}, the torsion restricted to $H$ has the form
\begin{equation}  \label{torha}
T(X,Y)=-[X,Y]_{|V}=2\sum_{s=1}^3\omega_s(X,Y)\xi_s, \qquad X,Y\in H.
\end{equation}

The properties of the Biquard connection are encoded in the properties of
the torsion endomorphism $T_{\xi}=T(\xi,.) : H\rightarrow H, \quad \xi\in V$%
. Decomposing the endomorphism $T_{\xi}\in(sp(n)+sp(1))^{\perp}$ into its
symmetric part $T^0_{\xi}$ and skew-symmetric part $b_{\xi}$, $%
T_{\xi}=T^0_{\xi} + b_{\xi} $, we summarize the description of the torsion
due to O. Biquard in the following Proposition. 

\begin{prop}
\cite{Biq1}\label{torb} The torsion $T_{\xi}$ is completely trace-free,
\begin{equation*}  \label{torb0}
tr\, T_{\xi}=g\,(\,T_{\xi}(e_a),e_a)=0, \quad tr\, T_{\xi}\circ
I=g\,(\,T_{\xi}(e_a),Ie_a)=0, \quad I\in Q,
\end{equation*}
where $e_1\dots e_{4n}$ is an orthonormal basis of $H$. Decomposing the
torsion into symmetric and antisymmetric parts, $T_{\xi_i}=T^0_{\xi_i}+b_{%
\xi_i}, \quad i=1,2,3$, we have: the symmetric part of the torsion has the
properties
\begin{gather*}  \label{tors1}
T^0_{\xi_i}I_i=-I_iT^0_{\xi_i} \\
\begin{aligned}\label{tors2}
I_2(T^0_{\xi_2})^{+--}=I_1(T^0_{\xi_1})^{-+-},\quad
I_3(T^0_{\xi_3})^{-+-}=I_2(T^0_{\xi_2})^{--+},\quad
I_1(T^0_{\xi_1})^{--+}=I_3(T^0_{\xi_3})^{+--}; \end{aligned}
\end{gather*}
the skew-symmetric part can be represented in the following way
\begin{equation*}  \label{toras}
b_{\xi_i}=I_iu,
\end{equation*}
where $u$ is a traceless symmetric (1,1)-tensor on $H$ which commutes with $%
I_1,I_2,I_3$.

If $n=1$ then the tensor $u$ vanishes identically, $u=0$ and the torsion is
a symmetric tensor, $T_{\xi}=T^0_{\xi}$.
\end{prop}

The covariant derivative of the quaternionic contact structure with respect
to the Biquard connection and the covariant derivative of the distribution $%
V $ are given by
\begin{equation}  \label{der}
\nabla I_i=-\alpha_j\otimes I_k+\alpha_k\otimes I_j,\qquad
\nabla\xi_i=-\alpha_j\otimes\xi_k+\alpha_k\otimes\xi_j,
\end{equation}
where the $sp(1)$-connection 1-forms $\alpha_s$ on $H$ are given by \cite%
{Biq1}
\begin{gather}  \label{coneforms}
\alpha_i(X)=d\eta_k(\xi_j,X)=-d\eta_j(\xi_k,X), \quad X\in H, \quad \xi_i\in
V,
\end{gather}
while the $sp(1)$-connection 1-forms $\alpha_s$ on the vertical space $V$
are calculated in \cite{IMV}
\begin{multline}  \label{coneform1}
\alpha_i(\xi_s)\ =\ d\eta_s(\xi_j,\xi_k) -\ \delta_{is}\left(\frac{Scal}{
16n(n+2)}\ +\ \frac12\,\left(\, d\eta_1(\xi_2,\xi_3)\ +\
d\eta_2(\xi_3,\xi_1)\ + \ d\eta_3(\xi_1,\xi_2)\right)\right),
\end{multline}
where $s\in\{1,2,3\}$. {The vanishing of the $sp(1)$-connection 1-forms on $%
H $ is equivalent to the vanishing of the torsion endomorphism of the
Biquard connection, see \cite{IMV}.}

\subsection{The qc-Einstein condition and Bianchi identities}

We explain briefly the consequences of the Bianchi identities and the notion
of qc-Einstein manifold introduced in \cite{IMV} since it plays a crucial
role in solving the Yamabe equation in the quaternionic seven dimensional
sphere. For more details see \cite{IMV}.

Let $R=[\nabla,\nabla]-\nabla_{[\ ,\ ]}$ be the curvature tensor of $\nabla$%
. The Ricci tensor and the scalar curvature $Scal$ of the Biquard
connection, called \emph{qc-Ricci tensor} and \emph{qc-scalar curvature},
respectively, are defined by
\begin{equation*}  \label{e:horizontal ricci}
Ric(X,Y)={g(R(e_a,X)Y,e_a)}, \quad X,Y \in H, \qquad
Scal=Ric(e_a,e_a)=g(R(e_b,e_a)e_a,e_b).
\end{equation*}
According to \cite{Biq1} the Ricci tensor restricted to $H$ is a symmetric
tensor. If the trace-free part of the qc-Ricci tensor is zero we call the
quaternionic structure \emph{a qc-Einstein manifold} \cite{IMV}. It is shown
in \cite{IMV} that the qc-Ricci tensor is completely determined by the
components of the torsion. First, recall the notion of the $Sp(n)Sp(1)$%
-invariant trace-free symmetric 2-tensors $T^0, U$ on $H$ introduced in \cite%
{IMV} by
\begin{gather*}  \label{tor}
T^0(X,Y)\overset{def}{=}g((T_{\xi_1}^{0}I_1+T_{\xi_2}^{0}I_2+T_{%
\xi_3}^{0}I_3)X,Y),\quad U(X,Y)\overset{def}{=}g(uX,Y), \quad X,Y\in H.
\end{gather*}
The tensor $T^0$ belongs to $[-1]$-eigenspace while $U$ is in the
[3]-eigenspace of the operator $\dag$ {given by \eqref{e:cross}}, i.e., they
have the properties:
\begin{gather} \label{propt}
T^0(X,Y)+T^0(I_1X,I_1Y)+T^0(I_2X,I_2Y)+T^0(I_3X,I_3Y)=0, \\\label{propu}
3U(X,Y)-U(I_1X,I_1Y)-U(I_2X,I_2Y)-U(I_3X,I_3Y)=0.
\end{gather}
\noindent{\ Theorem~1.3, Theorem~3.12 and Corollary~3.14 in \cite{IMV} imply:%
} 

\begin{thrm}
\cite{IMV}\label{sixtyseven} {Let $(M^{4n+3},g,\mathbb{Q})$ be a
quaternionic contact} $(4n+3)$-dimensional manifold, $n>1$ . For any $X,Y
\in H$ the qc-Ricci tensor and the qc-scalar curvature satisfy
\begin{equation*}  \label{sixtyfour}
\begin{aligned} Ric(X,Y) \ & =\ (2n+2)T^0(X,Y)
+(4n+10)U(X,Y)+\frac{Scal}{4n}g(X,Y)\\ Scal\ & =\
-8n(n+2)g(T(\xi_1,\xi_2),\xi_3) \end{aligned}
\end{equation*}
For $n=1$ the above formulas hold with $U=0$.

In particular, the qc-Einstein condition is equivalent to the vanishing of
the torsion endomorphism of the Biquard connection. If $Scal\not=0$ {\ the
latter} holds exactly when the qc-structure is 3-Sasakian up to a
multiplication by a constant and an $SO(3)$-matrix with smooth entries.
\end{thrm}

{For the last part of the above Theorem, we remind that a (4n+3)-dimensional
Riemannian manifold $(M,g)$ is called 3-Sasakian if the cone metric $%
g_N=t^2g+dt^2$ on $N=M\times \mathbb{R}^+$ is a hyperk\"ahler metric,
namely, it has holonomy contained in $Sp(n+1)$. }

The Ricci 2-forms $\rho_s, \ s=1,2,3$ of a quaternionic contact structure
are defined by
\begin{equation*}
4n\ \rho_s(B,C)=g(R(B,C)e_a,I_se_a),\ B,C\in \Gamma(TM).
\end{equation*}
\noindent {\ For ease of reference, in the following Theorem we summarize
the properties of the Ricci 2-forms, the scalar curvature and the torsion
evaluated on the vertical space established in Lemma~3.11, Corollary~3.14
Proposition~4.3 and Proposition~4.4 of \cite{IMV} }. 

\begin{thrm}
\cite{IMV}\label{summ} The Ricci 2-forms satisfy
\begin{align}
& \rho_1(X,Y)\ =\ 2g((T_{\xi_2}^0)^{--+}I_3X,Y)-2g(I_1uX,Y)-\frac{Scal}{%
8n(n+2)}\omega_1(X,Y),  \notag \\
& \rho_2(X,Y)\ =\ 2g((T_{\xi_3}^0)^{+--}I_1X,Y)-2g(I_2uX,Y)-\frac{Scal}{%
8n(n+2)}\omega_2(X,Y),  \label{rhoh} \\
& \rho_3(X,Y)\ =\ 2g((T_{\xi_1}^0)^{-+-}I_2X,Y)-2g(I_3uX,Y)-\frac{Scal}{
8n(n+2)}\omega_3(X,Y).  \notag
\end{align}
\begin{align}
& \rho_i(X,\xi_i)=-\frac{X(Scal)}{32n(n+2)}+
\frac12(\omega_i([\xi_j,\xi_k],X)-\omega_j([\xi_k,\xi_i],X)-\omega_k([\xi_i,%
\xi_j],X)),  \notag \\
& \rho_i(X,\xi_j)=\omega_j([\xi_j,\xi_k],X), \quad
\rho_i(X,\xi_k)=\omega_k([\xi_j,\xi_k],X),  \label{sixtyeight} \\
&
\rho_i(I_kX,\xi_j)=-\rho_i(I_jX,\xi_k)=g(T(\xi_j,\xi_k),I_iX)=\omega_i([%
\xi_j,\xi_k],X),  \notag \\
&\rho_i(\xi_i,\xi_j)+\rho_k(\xi_k,\xi_j)=\frac{1}{8n(n+2)}\xi_j(Scal).
\label{ricvert}
\end{align}
The torsion of the Biquard connection restricted to $V$ satisfies the
equality
\begin{equation}  \label{cornv}
T(\xi_i,\xi_j)=-\frac{Scal}{8n(n+2)}\xi_k-[\xi_i,\xi_j]_H,
\end{equation}
{\ where $[\xi_i,\xi_j]_H$ denotes the projection on $H$ parallel to the
vertical space $V$.}
\end{thrm}

We also recall the definition of the $Sp(n)Sp(1)$-invariant vector field ${A}
$, which appeared naturally in the Bianchi identities investigated in \cite%
{IMV}
\begin{equation*}  \label{d:A}
{A} \ = \ I_1[\xi_2,\xi_3]+I_2[\xi_3,\xi_1]+I_3[\xi_1,\xi_2].
\end{equation*}
\noindent We shall denote with the same letter the corresponding horizontal
one-form, i.e.,
\begin{equation*}
{A} (X)\ = \ g(I_1[\xi_2,\xi_3]+I_2[\xi_3,\xi_1]+I_3[\xi_1,\xi_2], X).
\end{equation*}
The horizontal divergence $\nabla^*P$ of a (0,2)-tensor field $P$ on $M$
with respect to Biquard connection is defined to be the (0,1)-tensor field
\begin{equation*}
\nabla^*P(.)=(\nabla_{e_a}P)(e_a,.).
\end{equation*}
Then we conclude from \cite[Theorem 4.8]{IMV}, that

\begin{thrm}
\cite{IMV}\label{bidiver} On a $(4n+3)$-dimensional QC manifold with
constant qc-scalar curvature we have the formulas
\begin{equation}  \label{div:To}
\nabla^*T^0=(n+2){A}, \qquad \nabla^*U=\frac{1-n}{2}{A}
\end{equation}
\end{thrm}


\section{Conformal transformations}

\label{s:conf transf}

Note that a conformal quaternionic contact transformation between two
quaternionic contact manifold is a diffeomorphism $\Phi$ which satisfies
\begin{equation*}
\Phi^*\eta=\mu\ \Psi\cdot\eta,
\end{equation*}
for some positive smooth function $\mu$ and some matrix $\Psi\in SO(3)$ with
smooth functions as entries and $\eta$ is an $\mathbb{R}^3$-valued one form,
$\eta=(\eta_1,\eta_2,\eta_3)^t$ is a column vector with entries one-forms.
The Biquard connection does not change under rotations, i.e., the Biquard
connection of $\Psi\cdot\eta$ and $\eta$ coincide. Hence, studying conformal
transformations we may consider only transformations $\Phi^*\eta\ =\ \mu\
\eta$.

Let $h$ be a positive smooth function on a QC manifold $(M, \eta)$. Let $%
\bar\eta=\frac{1}{2h}\eta$ be a conformal deformation of the QC structure $%
\eta$. We will denote the objects related to $\bar\eta$ by over-lining the
same object corresponding to $\eta$. Thus, $d\bar\eta=-\frac{1}{2h^2}%
\,dh\wedge\eta\ +\ \frac{1}{2h\,}d\eta$ and $\bar g=\frac{1}{2h}g$. The new
triple $\{\bar\xi_1,\bar\xi_2,\bar\xi_3\}$ is determined by the conditions
defining the Reeb vector fields. We have
\begin{equation}  \label{New19}
\bar\xi_s\ =\ 2h\,\xi_s\ +\ I_s\nabla h,\quad s=1,2,3,
\end{equation}
where $\nabla h$ is the horizontal gradient defined by $g(\nabla
h,X)=dh(X),\quad X\in H$.

\noindent The components of the torsion tensor transform according to the
following formulas from \cite[Section 5]{IMV}
\begin{gather}  \label{e:T^0 conf change}
\overline T^0(X,Y) \ =\ T^0(X,Y)\ +\ h^{-1}\,[\nabla dh]_{[sym][-1]}(X,Y),
\\\label{e:U conf change} \bar U(X,Y) \ =\ U(X,Y)\ +\ (2h)^{-1}[\ {\nabla
dh}-2h^{-1}dh\otimes dh]_{[3][0]}(X,Y),
\end{gather}
where the symmetric part is given by, {\ cf. \eqref{comm-rel},}
\begin{equation*}  \label{symdh}
[\ {\nabla dh}]_{[sym]}(X,Y)\ =\ \ {\nabla dh}(X,Y)\ + \ \sum_{s=1}^3
dh(\xi_s)\,\omega_s(X,Y)
\end{equation*}
{and $_{[3][0]}$ indicates the trace free part of the [3]-component of the
corresponding tensor. }In addition, the qc-scalar curvature changes
according to the formula \cite{Biq1}
\begin{equation}  \label{e:conf change scalar curv}
\overline {\text{Scal}}\ =\ 2h\,(\text{Scal})\ -\ 8(n+2)^2\,h^{-1}|\nabla
h|^2\ +\ 8(n+2)\,\triangle h.
\end{equation}
The following vectors will be important for our considerations,
\begin{equation}  \label{d:A_s}
A_i\ =\ I_i[\xi_j, \xi_k],\qquad\text{hence}\quad A\ = \ A_1\ +\ A_2\ +\ A_3.
\end{equation}

\begin{lemma}
\label{l:transf of A} Let $h$ be a positive smooth function on a QC manifold
$(M, g, \mathbb{Q})$ {\ with constant qc-scalar curvature $Scal=16n(n+2)$}
and $\bar\eta\ =\ \frac{1}{2h}\, \eta$ a conformal deformation of the qc
structure $\eta$. {If $\bar\eta$ is a 3-Sasakian structure, then} we have
the formulas
\begin{multline}  \label{e:A_s}
A_1(X)\ =\ -\frac12 h^{-2}dh(X)\ -\ \frac 12h^{-3}\lvert \nabla h
\rvert^2dh(X) \\
-\ \frac 12 h^{-1}\Bigl (\ {\nabla dh} (I_2X, \xi_2)\ +\ \ {\nabla dh}
(I_3X, \xi_3) \Bigr ) +\ \frac 12 h^{-2}\Bigl (dh(\xi_2)\,dh (I_2X)\ +\
dh(\xi_3)\,dh (I_3X) \Bigr ) \\
+\ \frac 14 h^{-2}\Bigl ( \ {\nabla dh} (I_2X, I_2 \nabla h)\ +\ \ {\nabla dh%
} (I_3X, I_3 \nabla h) \Bigr ).
\end{multline}
The expressions for $A_2$ and $A_3$ can be obtained from the above formula
by a cyclic permutation of $(1,2,3)$. Thus, we have also
\begin{multline*}  \label{e:A}
A(X)\ =\ -\frac32 h^{-2}dh(X)\ -\ \frac 32h^{-3}\lvert \nabla h \rvert^2dh(X)
\\
-\ h^{-1}\sum_{s=1}^3\ {\nabla dh} (I_sX, \xi_s)\ +\
h^{-2}\sum_{s=1}^3dh(\xi_s)\,dh (I_sX)\ +\ \frac 12 h^{-2}\sum_{s=1}^3\ {%
\nabla dh} (I_sX, I_s \nabla h)\
\end{multline*}
\end{lemma}

\begin{proof}
First we calculate the $sp\,(1)$-connection 1-forms of the Biquard
connection $\nabla$. For a 3-Sasaki structure we have $d\bar\eta_i(\bar%
\xi_j,\bar\xi_k)=2, \quad \bar\xi_i\lrcorner d\bar\eta_i=0$, the non-zero $%
sp(1)$-connection 1-forms are $\bar\alpha_i(\bar\xi_i)=-2, i=1,2,3$, and the
qc-scalar curvature $\overline {Scal}=16n(n+2)$ (see [Example 4.12,\cite{IMV}%
]). Then \eqref{New19}, \eqref{coneforms}, and \eqref{coneform1} yield
\begin{equation}  \label{con1forms}
\begin{aligned} 2d\eta_i(\xi_j,\xi_k)\ =\ 2h^{-1}\ +\ h^{-2}\lVert
dh\rVert^2, \qquad\qquad \alpha_i(X)\ =\ -h^{-1}dh(I_iX),\hskip.4truein\\
\alpha_i(\xi_j)\ =\ -h^{-1}dh(\xi_k)\ =\ -\alpha_j(\xi_i),\qquad
4\alpha_i(\xi_i)\ =\ -4\ -\ 2h^{-1}\ -\ h^{-2}\lVert dh\rVert^2.
\end{aligned}
\end{equation}
\noindent 
From the 3-Sasakian assumption the commutators are $[\bar\xi_i,\bar\xi_j]\
=\ -2\bar\xi_k$. Thus, for $X\in H$ taking also into account \eqref{New19}
we have
\begin{equation*}
g([\bar\xi_1,\bar\xi_2],\, I_3X)\ =\ -2\,g(\bar\xi_3, \, I_3X)\ =\ -2\,
g(2\,h\,\xi_3 \ +\ I_3\,\nabla h, \, I_3\,X)\ =\ -2\,dh\, (X).
\end{equation*}
\noindent Therefore, using again \eqref{New19}, we obtain
\begin{multline}  \label{e:com1}
-2\,dh(X)\ =\ g([\bar\xi_1,\, \bar\xi_2],\, I_3X)\ =\ g\bigl(\,[2h\xi_1\ +\
I_1\,\nabla h,2h\xi_2\ +\ I_2\,\nabla h],\, I_3X\bigr) \\
= \ -4h^2\,A_3(X)\ +\ 2hg(\,[\xi_1,\, I_2\nabla h],\, I_3\,X)\ +\ 2hg(\,
[I_1\nabla h,\xi_2],\, I_3\,X) \\
+\ g(\,[I_1\,\nabla h,\, I_2\,\nabla h],I_3\,X).
\end{multline}
\noindent The last three terms are transformed as follows. The first equals
\begin{multline*}
g\,(\,[\xi_1,\, I_2\nabla h],I_3\,X)\ = \ g\,(\, (\nabla_{\xi_1}\, I_2)\,
\nabla h\ +\ I_2\, \nabla_{\xi_1}\, \nabla h, \, I_3X) \ -\ g\,(\,T(\xi_1,
\, I_2 \nabla h), \, I_3X) \\
=\ -\alpha_3(\xi_1)\, dh(I_2X)\ +\ \alpha_1(\xi_1)\, dh(X)\ -\ \nabla dh
\,(\xi_1, I_1X)\ -\ g\,(\,T(\xi_1, \, I_2 \nabla h), \, I_3X),
\end{multline*}
\noindent {where we use \eqref{der} and the fact that $\nabla$ preserves the
splitting $H\oplus V$.} The second term is
\begin{multline*}
g\,( \,[I_1\nabla h,\xi_2],I_3\,X)\ =\ \alpha_2(\xi_2)\, dh(X)\ +\
\alpha_3(\xi_2)\, dh(I_1X) \ -\ \nabla dh\, (\xi_2, I_2X) \\
-\ g\,(\,T( I_1 \nabla h,\, \xi_2), \, I_3X),
\end{multline*}
\noindent and finally
\begin{multline*}
g\,(\, [I_1\nabla h,I_2\nabla h],I_3\,X)\ =\ -\alpha_3(I_1 \nabla
h)\,dh(I_2X)\ +\ \alpha_1(I_1 \nabla h)\,dh(X)\ -\ \nabla dh\,(I_1\nabla
h,I_1X) \\
+\ \alpha_2(I_2 \nabla h)\,dh(X)\ +\ \alpha_3(I_2 \nabla h)\,dh(I_1X)\ -\ \nabla
dh\,(I_2\nabla h,I_2X).
\end{multline*}
\noindent Next we apply \eqref{con1forms} to the last three equalities, then
substitute their sum into \eqref{e:com1}, after which we use the commutation
relations
\begin{equation}  \label{comm-rel}
\begin{aligned} \nabla dh\,(X,Y)-\nabla dh\,(Y,X)\ =\ -dh(T(X,Y))=\
-2\sum_{s=1}^3\omega_s(X,Y)\,dh (\xi_s), \\ \nabla dh\,(X,\xi)-\nabla
dh\,(\xi,X)\ =\ -dh(T(X,\xi)), \qquad X, \, Y \in H,\quad \xi\in V.
\end{aligned}
\end{equation}
The result is the following identity
\begin{multline}  \label{e:A_3 comp}
4h^2A_3(X)\ =\ (-4h\ +\ h^{-1}\lVert \nabla h\rVert^2)\, dh(X) \\
-2h\,\bigl[\, \nabla dh\, (I_1X,\xi_1)\ +\ \nabla dh\, (I_2X,\xi_2)\bigr]-\ %
\bigl [\, \nabla dh\, (I_1X,I_1\nabla h))\ +\ \nabla dh\, (I_2X,I_2\nabla h)%
\bigr] \\
+\ 2\, \bigl[\, dh(\xi_1)\,dh(I_1X)\ +\ dh(\xi_2)\,dh(I_2X)\ +\
2\,dh(\xi_3)\,dh(I_3X)\bigr] \\
+\ 2h\,\bigl[ \, T(\xi_1,I_1X, \nabla h)\ +\ T(\xi_2,I_2X, \nabla h)\ -\
T(\xi_1,I_2 X, I_3\nabla h)\ +\ T(\xi_2,I_1X, I_3 \nabla h)\bigr ],
\end{multline}
{where $T(\xi, X, Y)=g(T_\xi X, Y)$ for a vertical vector $\xi$ and
horizontal vectors $X$ and $Y$. With the help of Proposition~\ref{torb} we
decompose the torsions into symmetric and anti-symmetric part $T_{\xi_i}=
T^0_{\xi_i} + I_iU$, $i=1,\ 2,\ 3$, and then express the symmetric parts of
the torsion terms in the form $T^0_{\xi_1}=(T^0_{\xi_1})^{--+}+(T^0_{%
\xi_1})^{-+-}$, $T^0_{\xi_2}=(T^0_{\xi_2})^{--+}+(T^0_{\xi_2})^{+--}$.
Hence, using $T^{0^{--+}}=2(T^0_{\xi_2})^{+--}I_2=2(T^0_{\xi_1})^{-+-}I_1$
etc., which follows again from Proposition~\ref{torb}, the sum of the
torsion terms in \eqref{e:A_3 comp} can be seen to equal $2T^{0^{--+}}(X ,\,
\nabla h)\ -\ 4\,U(X, \nabla h)$. } This allows us to rewrite
\eqref{e:A_3
comp} in the form
\begin{multline}  \label{e:A_3 comp 2}
4A_3(X)\ =\ (-4h^{-1}\ +\ h^{-3}\lVert \nabla h\rVert^2)\,dh(X) -\
2\,h^{-1}\,\bigl[\, \nabla dh\, (I_1X,\xi_1)\ +\ \nabla dh\, (I_2X,\xi_2)%
\bigr] \\
+\ 2\,h^{-2}\,\bigl[\, dh(\xi_1)\,dh(I_1X)\ +\,dh(\xi_2)\,dh(I_2X)\ +\
2\,dh(\xi_3)\,dh(I_3X)\bigr] \\
-\ h^{-2}\,\bigl[\,\nabla dh\, (I_1X,I_1\nabla h)\ +\ \nabla dh\,
(I_2X,I_2\nabla h)\bigr]\ +\ 4h^{-1}\,\bigl[\,(T^{0^{--+}}(\nabla h,X)\ -\
2U(\nabla h,X)\bigr].
\end{multline}
\noindent Using \eqref{e:T^0 conf change} the $T^{0^{--+}} $ component of
the torsion can be expressed by $h$ as follows, {see \eqref{New21} and
\eqref{propt},}
\begin{multline*}
4T^{0^{--+}}(\nabla h,X)\ =\ T^0(\nabla h,X)\ -\ T^0(I_1\nabla h,I_1X)\ -\
T^0(I_2\nabla h,I_2X)\ +\ T^0(I_3\nabla h,I_3X) \\
= \ -h^{-1}\Bigl \{ \lbrack \nabla dh]_{[-1]} (\nabla h , X) - [\nabla
dh]_{[-1]} (I_1\nabla h,I_1X) -[\nabla dh]_{[-1]}(I_2\nabla h,I_2X)+ [\nabla
dh]_{[-1]}(I_3\nabla h,I_3X)\Bigr \} \\
-h^{-1}\, \sum_{s=1}^3\Bigl \{ dh(\xi_s)\,\Bigl [ g(I_s \nabla h,X)\ -\
g(I_sI_1\nabla h,I_1X) -\ g(I_sI_2\nabla h,I_2X)\ +\ g(I_sI_3\nabla h,I_3X) %
\Bigr ] \Bigr \} \\
=\ -h^{-1}\, \Bigl \{ \nabla dh\,(\nabla h , X)\ -\ \nabla dh\,(I_1\nabla h
, I_1X)\ -\ \nabla dh\,(I_2\nabla h , I_2X)\ +\ \nabla dh\,(I_3\nabla h ,
I_3X) \Bigr \} \\
+\ 4h^{-1}\, dh(\xi_3)\, dh(I_3X).
\end{multline*}
{\ Invoking equation \eqref{comm-rel} we can put $\nabla h$ in second place
in the Hessian terms, thus, proving the formula}
\begin{multline}  \label{e:T--+}
4T^{0^{--+}}(\nabla h,X)\ =\ -\ 4h^{-1}\, dh(\xi_3)\, dh(I_3X) \\
-h^{-1}\, \Bigl \{ \nabla dh\,( X, \nabla h )\ -\ \nabla dh\,( I_1X,
I_1\nabla h )\ -\ \nabla dh\,( I_2X, I_2\nabla h )\ +\ \nabla dh\,( I_3X,
I_3\nabla h ) \Bigr \} .
\end{multline}
{On the other hand, \eqref{propu}, \eqref{e:U conf change} and the Yamabe
equation } \eqref{e:conf change scalar curv} give
\begin{multline}  \label{e:U by h}
8\,U(\nabla h, X)\ =\ -h^{-1}\, \Bigl \{ \nabla dh\,(\nabla h , X)\ + \
\sum_{s=1}^3 \nabla dh\,(I_s\nabla h , I_sX) \\
\hskip1.9truein - 2h^{-1}\lVert \nabla h\rVert^2\, dh (X)\ -\ \frac {%
\triangle h}{n}\, dh(X)\ +\ 2h^{-1} \frac {\lVert \nabla h\rVert^2}{n}\, dh
(X) \Bigr \} \\
=\ -h^{-1}\, \Bigl \{ \nabla dh\,(\nabla h , X)\ + \ \sum_{s=1}^3 \nabla
dh\,(I_s\nabla h , I_sX) \Bigr \} \\
-\ h^{-1}\, \Bigl \{ - 2 h^{-1}\lVert \nabla h\rVert^2\, dh (X)\ -\ \frac {%
2n-4nh+ (n+2)h^{-1}\lVert \nabla h\rVert^2}{n}\, dh(X)\ +\ 2 h^{-1} \frac {%
\lVert \nabla h\rVert^2}{n}\, dh (X) \Bigr \} \\
=\ -h^{-1}\, \Bigl \{ \nabla dh\,(X, \nabla h)\ + \ \sum_{s=1}^3 \nabla
dh\,( I_sX, I_s\nabla h ) \Bigr \} -\ h^{-1}\, \bigl (-3h^{-1}\lVert \nabla
h\rVert^2\ -\ 2\ +\ 4h \bigr )dh (X).
\end{multline}
\noindent Substituting the last two formulas in \eqref{e:A_3 comp 2} gives {%
\ $A_3$ in the form} of \eqref{e:A_s} written for $A_1$, cf. the paragraph
after \eqref{e:A_s}.
\end{proof}


\section{Divergence formulas}

\label{s:div formulas} We shall need the divergences of various vector/forms
through the almost complex structures, so we start with a general formula
valid for any horizontal vector/form $A$. Let $\{e_1, \dots, e_{4n}\}$ be an
orthonormal basis of $H$. The divergence of $I_1 A$ is
\begin{equation*}
\nabla^* (I_1 A)\ \equiv\ (\nabla_{e_a} (I_1A))(e_a)\ =\ - (\nabla_{e_a}
A)(I_1e_a)\ -\ A((\nabla_{e_a}I_1)e_a),
\end{equation*}
recalling $I_1A(X) = -A(I_1X)$.

{\ We say that an orthonormal frame\newline
\centerline{$\{e_1,e_2=I_1e_1,e_3=I_2e_1,e_4=I_3e_1,\dots,
e_{4n}=I_3e_{4n-3}, \xi_1, \xi_2, \xi_3 \}$}\newline is a
qc-normal frame (at a point) if the connection 1-forms of the
Biquard connection vanish (at that point). Lemma~4.5 in \cite{IMV}
asserts that a qc-normal frame exists at each point of a QC
manifold.} With respect to a qc-normal frame the above divergence
reduces to
\begin{equation*}
\nabla^* (I_1 A)\ \ =\ - (\nabla_{e_a} A)(I_1e_a).
\end{equation*}

\begin{lemma}
\label{l:div of I_sA} Suppose $(M, \eta, \mathbb{Q})$ is a quaternionic
contact manifold with constant qc-scalar curvature. For any function $h$ we
have the following formulas
\begin{equation*}
\begin{aligned} \nabla^*\, \Bigl (\sum_{s=1}^3 dh(\xi_s) I_sA_s\Bigr )\ =\
\sum_{s=1}^3 \ \nabla dh\,(I_s e_a, \xi_s)A_s(e_a)\\ \nabla^*\, \Bigl
(\sum_{s=1}^3 dh(\xi_s) I_sA \Bigr )\ =\ \sum_{s=1}^3 \ \nabla dh\,(I_s e_a,
\xi_s)A(e_a). \end{aligned}
\end{equation*}
\end{lemma}

\begin{proof}
Using the identification of the 3-dimensional vector spaces
spanned by $\{\xi_1,\xi_2, \xi_3\} $ and $\{I_1,I_2,I_3 \}$ with
$\mathbb{R}^3$, the restriction of the action of $Sp(n)Sp(1)$ to
this spaces can be identified with  the action of the group
$SO(3)$.  With this in mind, one verifies easily that the 1-forms
$A$, $\sum_{s=1}^3 dh(\xi_s) I_sA_s=-\sum_{i=1}^3 dh(\xi_i)
[\xi_j,\xi_k]$ and $\sum_{s=1}^3 dh(\xi_s) I_sA$ are $Sp(n)Sp(1)$
invariant on $\mathbb H$. Thus, it is sufficient to compute their
divergences in a qc-normal frame. To avoid the introduction of new
variables, in this proof, we shall assume that $\{e_1, \dots,
e_{4n}, \xi_1, \xi_2, \xi_3 \}$ is a qc-normal frame.

{\ We apply \eqref{cornv}.} Using that the Biquard connection preserves the
splitting of $TM$, we find
\begin{multline*}
\nabla^* [\xi_1, \xi_2]\ =\ -g( \nabla_{e_a}\bigl ( T(\xi_1, \xi_2)\bigr ),
e_a ) \\
=\ -g(\,(\nabla_{e_a} T)\, (\xi_1, \xi_2), e_a )\ -\ g( \,
T(\nabla_{e_a}\xi_1, \xi_2), e_a )\ -\ g(\, T(\xi_1, \nabla_{e_a}\xi_2), e_a
).
\end{multline*}
From Bianchi's identity we have {\ ($\sigma_{A,B,C}$ means a cyclic sum
over $(A,B,C)$)} 
\begin{multline*}
g((\nabla_{e_a} T) (\xi_1, \xi_2), e_a )\ =\ - g((\nabla_{\xi_1} T)(\xi_2,
e_a), e_a)\ -\ g((\nabla_{\xi_2} T) (e_a, \xi_1), e_a) \\
-\ g(\sigma_{e_a, \xi_1, \xi_2} \Bigl \{ T(T(e_a, \xi_1),\xi_2)\Bigr \},e_a)
\ + g( \sigma_{e_a, \xi_1, \xi_2} \Bigl \{R(e_a, \xi_1)\xi_2\Bigr \}, e_a) \\
=\ - g(T(\, T(e_a, \xi_1), \xi_2), e_a)\ -\ g(T (\,T(\xi_1, \xi_2), e_a),
e_a)\ -\ g(T(T(\xi_2, e_a), \xi_1), e_a) \\
=\ g(\,T (T({\xi_1},e_a), \xi_2), e_a)\ -\ g(\,T (T({\xi_2},e_a), \xi_1),
e_a)\ -\ g(\,T (T(\xi_1, \xi_2), e_a), e_a),
\end{multline*}
{taking into account that as mappings on $H$ the torsion tensors $T({\xi_i}%
,X)$ and the curvature tensor $R(\xi_1, \xi_2)$ are traceless}, {so $%
g((\nabla_{\xi_1} T)(\xi_2, e_a), e_a)$ and $g(R(\xi_1, \xi_2) e_a, e_a)=0$,
while the connection preserves the splitting,} to obtain the next to last
line. The last term is equal to zero as
\begin{equation*}
g(T (T(\xi_1, \xi_2), e_a), e_a)\ = \ g\bigl (T ( -\frac{\text{Scal}}{8n(n+2)%
}\xi_3\ -\ [\xi_1,\xi_2]_H, e_a), e_a\bigr )\ =\ -\frac{\text{Scal}}{8n(n+2)}%
\, g(T ( \xi_3, e_a), e_a)\ = 0,\qquad\qquad
\end{equation*}
taking into account that the torsion $T_{\xi_3}$ is traceless and $T (
[\xi_1,\xi_2]_H, e_a)$ is a vertical vector. On the other hand,
\begin{multline*}
g\,(\,T (T_{\xi_1}e_a, \xi_2), e_a)\ -\ g\,(\,T (T_{\xi_2}e_a, \xi_1), e_a)
\\
=\ -\ \Bigl[\,g(T(e_b,\xi_2),e_a)\, g(T(\xi_1,e_a),e_b)\ -\
g(T(e_b,\xi_1),e_a)\, g(T(\xi_2,e_a),e_b)\Bigr] \\
=\ \Bigl[g(T(\xi_2,e_b),e_a)\, g(T(\xi_1,e_a),e_b)\ - \
g(T(\xi_1,e_b),e_a)\, g(T(\xi_2,e_a),e_b)\Bigr]\ =\ 0.
\end{multline*}
\noindent The equalities $\nabla^*(I_1A_1)=\nabla^*(I_2A_2)=0$ with respect
to a qc-normal frame can be obtained similarly. Hence, the first formula in
Lemma~\ref{l:div of I_sA} follows.

We are left with proving the second divergence formula. Since the scalar
curvature is constant \eqref{sixtyeight} implies
\begin{equation}  \label{e:a1}
A(X)=-2\,\sum_{s=1}^3 \rho_s(X,\xi_s).
\end{equation}
\noindent Fix an $s\in\{1,2,3\}$. Working again in a qc-normal frame we have
\begin{equation*}
(\nabla_{e_a}A)(I_se_a)\ =\ -2\,\sum_{t=1}^3 (\nabla_{e_a} \rho_t)
(I_se_a,\xi_t).
\end{equation*}
\noindent A calculation involving {\ the expressions \eqref{rhoh}} and the
properties of the torsion shows that
\begin{equation}  \label{e:tr I rho}
\text{tr}\, (\rho_t\circ I_s)\ =\ -\frac {1}{2(n+2)}\, \delta_{st}\,\text{%
Scal}.
\end{equation}
\noindent The second Bianchi identity
\begin{multline*}
0\ =\ g((\nabla_{e_a} R)(I_se_a,\xi_t)e_b,I_te_b)\ +\ g((\nabla_{I_se_a}
R)(\xi_t,e_a)e_b,I_te_b)\ +\ g((\nabla_{\xi_t} R)(e_a,I_se_a)e_b,I_te_b) \\
+\ g(R(T(e_a,I_se_a),\xi_t)e_b,I_te_b)\ +\
g(R(T(I_se_a,\xi_t),e_a)e_b,I_te_b)\ +\ g(R(T(\xi_t,e_a),I_se_a)e_b,I_te_b),
\end{multline*}
\noindent together with the constancy of the qc-scalar curvature and %
\eqref{e:tr I rho} show that 
the third term on the right is zero and thus
\begin{equation*}
\sum_{t=1}^3\Bigl \{\, 2\,(\nabla_{e_a}\rho_t)(I_se_a,\xi_t)\ - \
2\,\rho_t(T(\xi_t,I_se_a),e_a)\ +\ \rho_t(T(e_a,I_se_a),\xi_t)\Bigr \}\ =\ 0.
\end{equation*}
\noindent Substituting \eqref{torha} in the above equality we come to the
equation
\begin{equation}  \label{bia1}
\sum_{t=1}^3(\nabla_{e_a}\rho_t)(I_se_a,\xi_t)\ =\
\sum_{t=1}^3\rho_t(T(\xi_t,I_se_a),e_a)-
4n\,\sum_{t=1}^3\rho_t(\xi_s,\xi_t)=0,
\end{equation}
\noindent where the vanishing of the second term follows from \eqref{ricvert}%
, while the vanishing of the first term is seen as follows. {\ Using the
standard inner product on $End(H)$}{\
\begin{equation*}
g(C,B)=tr(B^*C)=\sum_{a=1}^{4n}g(C(e_a),B(e_a)),
\end{equation*}
where $C,B\in End(H)$, $\{e_1,...,e_{4n}\}$ is a $g$-orthonormal basis of $H$%
,} the definition of $T^0_{\xi_s}$, {\ the formulas in Theorem~\ref{summ}
and Proposition~\ref{torb} imply}
\begin{multline*}
\sum_{s=1}^3\rho_s(T(\xi_s, I_1e_a), e_a) \\
=\ g(\rho_1, T^0_{\xi_1}I_1)\ +\ g(\rho_2, T^0_{\xi_2}I_1)\ +\ g(\rho_3,
T^0_{\xi_3}I_1)\ -\ g(\rho_1, u)\ -\ g(\rho_2, I_3u)\ +\ g(\rho_3, I_2u) \\
=\ g(\rho_1, T^0_{\xi_1}I_1)\ +\ g(\rho_2, T^0_{\xi_2}I_1)\ +\ g(\rho_3,
T^0_{\xi_3}I_1) \\
=\ g(2(T_{\xi_2}^0)^{--+}I_3\ -\ 2I_1u\ -\ \frac {Scal}{8n(n+2)}\,I_1,\,
T^0_{\xi_1}I_1) \\
+\ g(2(T_{\xi_3}^0)^{+--}I_1\ -\ 2I_2u\ -\ \frac {Scal}{8n(n+2)}\,I_2,\,
T^0_{\xi_2}I_1) +\ g(2(T_{\xi_1}^0)^{-+-}I_2\ -\ 2I_3u\ -\ \frac {Scal}{%
8n(n+2)}\,I_3,\, T^0_{\xi_3}I_1) \\
=\ -2g((T_{\xi_2}^0)^{--+}I_2, T^0_{\xi_1}) \ +\ 2 g((T_{\xi_3}^0)^{+--},
T^0_{\xi_2})\ +\ 2 g((T_{\xi_1}^0)^{-+-}I_3, T^0_{\xi_3}) \\
=\ 2 g((T_{\xi_3}^0)^{+--}, (T_{\xi_2}^0)^{+--})\ +\ 2
g((T_{\xi_1}^0)^{-+-}, I_3(T_{\xi_3}^0)^{+--}) \\
=\ 2 g(I_2(T_{\xi_3}^0)^{+--}, I_2(T_{\xi_2}^0)^{+--})\ -\ 2
g(I_1(T_{\xi_1}^0)^{-+-}, I_2(T_{\xi_3}^0)^{+--})\ =\ 0.
\end{multline*}
\noindent Renaming the almost complex structures shows that the same
conclusion is true when we replace $I_1$ with $I_2$ or $I_3$ in the above
calculation.

Finally, the second formula in Lemma~\ref{l:div of I_sA} follows from %
\eqref{e:a1} and \eqref{bia1}.
\end{proof}

We shall also need the following one-forms
\begin{equation}  \label{d:D_s}
\begin{aligned} D_1(X)\ =\ - h^{-1}{T^{0^{+--}}}(X,\nabla h)\\ D_2(X)\ =\ -
h^{-1}{T^{0^{-+-}}}(X,\nabla h)\\ D_3(X)\ =\ - h^{-1}{T^{0^{--+}}}(X,\nabla
h) \end{aligned}
\end{equation}
For simplicity, using the musical isomorphism, we will denote with $D_1, \,
D_2, \, D_3$ the corresponding (horizontal) vector fields, for example
\hspace{3mm} $g(D_1, X)=D_1(X) \qquad\forall X\in H.$ \hspace{3mm} Finally,
we set
\begin{equation}  \label{d:def of D}
D\ =\ D_1\ +\ D_2\ +\ D_3\ =\ - h^{-1}\,T^{0}(X,\nabla h).
\end{equation}

\begin{lemma}
\label{l:div of D} Suppose $(M, \eta)$ is a quaternionic contact manifold
with constant {\ qc-}scalar curvature $Scal=16n(n+2)$. {Suppose $%
\bar\eta=\frac1{2h}\eta$ has vanishing $[-1]$-torsion component $\overline
T^0=0$}. We have
\begin{equation*}
D(X)\ =\ \frac 14 h^{-2}\Bigl (3\ {\nabla dh} (X, \nabla h)\ -\
\sum_{s=1}^3\ {\nabla dh} (I_sX, I_s\nabla h) \Bigr )\ +\
h^{-2}\sum_{s=1}^3dh(\xi_s)\,dh (I_sX).
\end{equation*}
and the divergence of $D$ satisfies
\begin{equation*}
\nabla^*\, D\ =\ \lvert T^0 \rvert^2\ -h^{-1}g(dh,D)\ -\ h^{-1}
(n+2)\,g(dh,A).
\end{equation*}
\end{lemma}

\begin{proof}
a) The formula for $D$ follows immediately from \eqref{e:T^0 conf change}.

b) We work in a qc-normal frame. Since the scalar curvature is assumed to be
constant we use \eqref{div:To} to find
\begin{multline*}
\nabla^* D\ =\ -h^{-1}\, dh(e_a)D(e_a)\ -\ h^{-1}\nabla^*T^0(\nabla h)\ -\
h^{-1}T^0(e_a,e_b)\ {\nabla dh}(e_a,e_b) \\
=\ -h^{-1}\, dh(e_a)D(e_a)\ -\ h^{-1} (n+2)\,dh (e_a)A(e_a)\ -\ \,g(T^0,\,
h^{-1}\ {\nabla dh}) \\
=\ \lvert T^0 \rvert^2\ -h^{-1}dh(e_a)D(e_a)\ -\ h^{-1} (n+2)\,dh
(e_a)A(e_a),
\end{multline*}
using \eqref{e:T^0 conf change} in the last equality.
\end{proof}

Let us also consider the following one-forms (and corresponding vectors)
\begin{equation*}  \label{d:F_s}
F_s(X)\ =\ - h^{-1}\, {T^0}(X,I_s\nabla h), \quad X\in H\quad s=1,2,3.
\end{equation*}
\noindent From the definition of $F_1$ and \eqref{d:D_s} we find
\begin{multline*}
F_1(X)\ =\ - h^{-1}{T^0}(X,I_1\nabla h) \\
=\ - h^{-1}{T^0}^{+--}(X,I_1\nabla h)\ -\ h^{-1}{T^0}^{-+-}(X,I_1\nabla h)\
- h^{-1}{T^0}^{--+}(X,I_1\nabla h) \\
=\ h^{-1}{T^0}^{+--}(I_1X,\nabla h)\ -\ h^{-1}{T^0}^{-+-}(I_1X,\nabla h)\ -
h^{-1}{T^0}^{--+}(I_1X,\nabla h) \\
=\ -D_1(I_1X)\ +\ D_2(I_1 X)\ +\ D_3(I_1X).
\end{multline*}
\noindent Thus, the forms $F_s$ can be expressed by the forms $D_s$ as
follows
\begin{equation}  \label{e:F_s by D_s}
\begin{aligned} F_1(X)\ =\ -D_1(I_1X)\ +\ D_2(I_1 X)\ +\ D_3(I_1X)\\ F_2(X)\
=\ D_1(I_2X)\ -\ D_2(I_2 X)\ +\ D_3(I_2X)\\ F_3(X)\ =\ D_1(I_3X)\ +\ D_2(I_3
X)\ -\ D_3(I_3X). \end{aligned}
\end{equation}

\begin{lemma}
\label{p:div of F_s}Suppose $(M, \eta)$ is a quaternionic contact manifold
with constant {\ qc-}scalar curvature $\text{Scal}=16n(n+2)$. {\ Suppose $%
\bar\eta=\frac1{2h}\eta$ has vanishing $[-1]$-torsion component, $\overline
T^0=0$}. We have
\begin{multline*}
\nabla^*\, \Bigl (\sum_{s=1}^3 dh(\xi_s) F_s\Bigr )\ =\ \sum_{s=1}^3 \Bigl [%
\ \nabla dh\, (I_se_a,\xi_s)F_s(I_se_a)\Bigr] \\
+ \ h^{-1}\sum_{s=1}^3 \Bigl[dh(\xi_s)dh (I_se_a)D(e_a)\ +(n+2)\,dh(\xi_s)dh
(I_s e_a)\, A(e_a)\Bigr ].
\end{multline*}
\end{lemma}

\begin{proof}
We note that the vector $\sum_{s=1}^3 dh(\xi_s) F_s$ is an $Sp(n)Sp(1)$
invariant vector, hence, we may assume that $\{e_1, \dots, e_{4n}, \xi_1,
\xi_2, \xi_3 \}$ is a qc-normal frame. Since the scalar curvature is assumed
to be constant we can apply Theorem~\ref{bidiver}, thus $\nabla^* T^0\ =\
(n+2)A$. Turning to the divergence, we compute
\begin{multline}  \label{e:div Fs}
\nabla^* \Bigl ( \sum_{s=1}^3dh(\xi_s)F_s\Bigr )\ =\ \sum_{s=1}^3 \Bigl [\
\nabla dh\, (e_a,\xi_s)F_s(e_a)\Bigr ]\ -\ \sum_{s=1}^3
h^{-1}\,dh(\xi_s)\,\nabla^*T^0(I_s\nabla h) \\
+ \sum_{s=1}^3\Bigl [ \, h^{-2}\, dh(\xi_s)\, dh(e_a)T^0(e_a,I_se_b)\,dh
(e_b)\ -\ h^{-1}\, dh(\xi_s)\,T^0(e_a,I_se_b)\ \nabla dh\, (e_a,e_b)\Bigr ]
\\
=\ \sum_{s=1}^3 \Bigl [\ \nabla dh\, (e_a,\xi_s)F_s(e_a)\Bigr ]\ -\
\sum_{s=1}^3 h^{-1}\,dh(\xi_s)\,\nabla^*T^0(I_s\nabla h) \\
+\sum_{s=1}^3 \Bigl [\ h^{-1}\, dh(\xi_s)\,dh(I_se_a)\,D(e_a)\Bigr ] \\
=\ \sum_{s=1}^3 \Bigl [\ \nabla dh\, (e_a,\xi_s)F_s(e_a)\ + \
h^{-1}\,dh(\xi_s)\,dh (I_se_a)\,D(e_a)\ +\ h^{-1}(n+2)\,dh(\xi_s)\,dh (I_s
e_a)\, A(e_a)\Bigr ],
\end{multline}
\noindent using the symmetry of $T^0$ in the next to last
equality, and  the fact
$$T^0(e_a,I_1e_b)\,\
\nabla dh\, (e_a,e_b)\ =-\ h^{-1}\nabla
dh_{[sym][-1]}(e_a,I_1e_b)\Big[\nabla
dh_{[sym]}(e_a,e_b)-\sum_{s=1}^3dh(\xi_s)\omega_s(e_a,e_b)\Big]=0,$$
which is a consequence of \eqref{e:T^0 conf change}, the formula
for the symmetric part of $\nabla dh$ given after \eqref{e:U conf
change} and the zero traces of the [-1]-component. Switching to
the basis $\{I_se_a:a=1,\dots,4n\}$ in the first term of the
right-hand-side of \eqref{e:div Fs} completes the proof.
\end{proof}

{\ At this point we restrict our considerations to the 7-dimensional case,
i.e. $n=1$. Following is our main technical result.}{\ As mentioned in the
introduction, we were motivated to seek a divergence formula of this type
based on the Riemannian and CR cases of the considered problem. The main
difficulty was to find a suitable vector field with non-negative divergence
containing the norm of the torsion. The fulfilment of this task was
facilitated by the results of \cite{IMV}, which in particular showed that
similarly to the CR case, but unlike the Riemannian case, we were not able
to achieve a proof based purely on the Bianchi identities, see \cite[Theorem
4.8]{IMV}.} 

\begin{thrm}
\label{t:div formula} {\ Suppose $(M^7,\eta)$ is a quaternionic contact
structure conformal to a 3-Sasakian structure $(M^7,\bar\eta)$}, $%
\tilde\eta\ =\ \frac{1}{2h}\, \eta.$ If $Scal_{\eta}=Scal_{\tilde%
\eta}=16n(n+2)$, then with $f$ given by
\begin{equation*}
f\ = \ \frac 12\ +\ h\ +\ \frac 14 h^{-1}\lvert \nabla h \rvert^2,
\end{equation*}
\noindent the following identity holds
\begin{equation*}
\nabla^*\Bigl(fD\ +\ \sum_{s=1}^3 dh(\xi_s)\, F_s \ +\ 4\sum_{s=1}^3
dh(\xi_s)I_sA_s \ -\ \frac {10}{3}\sum_{s=1}^3 dh(\xi_s)\,I_s A \Bigr) =\
f\lvert T^0 \rvert^2\ +\ h\,\langle Q V,\, V\rangle .
\end{equation*}
Here, $Q$ is a positive semi-definite matrix and $V=( D_1, D_2, D_3,A_1,
A_2, A_3)$ with $A_s$, $D_s$ defined, correspondingly, in \eqref{d:A_s} and %
\eqref{d:D_s}.
\end{thrm}

\begin{proof}
{Using the formulas for the divergences of $D$, $\sum_{s=1}^{3}dh(\xi
_{s})\,F_{s}$, $\sum_{s=1}^{3}dh(\xi _{s})I_{s}A_{s}$ and $%
\sum_{s=1}^{3}dh(\xi _{s})\,I_{s}A$ given correspondingly in Lemmas \ref%
{l:div of D}, \ref{p:div of F_s} and \ref{l:div of I_sA} we have the
identity ($n=1$ here)
\begin{multline}
\nabla ^{\ast }\Bigl(\,fD\ +\ \sum_{s=1}^{3}dh(\xi _{s})\,F_{s}\ +\
4\sum_{s=1}^{3}dh(\xi _{s})I_{s}A_{s}\ -\ \frac{10}{3}\sum_{s=1}^{3}dh(\xi
_{s})\,I_{s}A\Bigr)  \label{e:diver} \\
=\ \Bigl (dh(e_{a})\ -\ \frac{1}{4}h^{-2}dh(e_{a})\lvert \nabla h\rvert
^{2}\ +\ \frac{1}{2}h^{-1}\,\ {\nabla dh}(e_{a},\nabla h)\Bigr )\,D(e_{a}) \\
+\ f\,\Bigl (\lvert T^{0}\rvert ^{2}\ -h^{-1}dh(e_{a})D(e_{a})\ -\
h^{-1}(n+2)\,dh(e_{a})A(e_{a})\Bigr ) \\
+\ \sum_{s=1}^{3}\ {\nabla dh}(I_{s}e_{a},\xi _{s})F_{s}(I_{s}e_{a})\ +\
h^{-1}\sum_{s=1}^{3}\Bigl [\,dh(\xi _{s})\,dh(I_{s}e_{a})\,D(e_{a})\ +\
(n+2)\,dh(\xi _{s})\,dh(I_{s}e_{a})\,A(e_{a})\Bigr ] \\
+\ 4\sum_{s=1}^{3}{\nabla dh}(I_{s}e_{a},\xi _{s})\,A_{s}(e_{a})\ -\frac{10}{%
3}\ \sum_{s=1}^{3}\ {\nabla dh}(I_{s}e_{a},\xi _{s})\,A(e_{a}) \\
=\ \Bigl (dh(e_{a})\ -\ \frac{1}{4}h^{-2}dh(e_{a})\lvert \nabla h\rvert
^{2}\ +\ \frac{1}{2}h^{-1}\,\ {\nabla dh}(e_{a},\nabla h)\Bigr )%
\,\sum_{t=1}^{3}D_{t}(e_{a}) \\
+\ f\,\Bigl (\lvert T^{0}\rvert ^{2}\ -h^{-1}dh(e_{a})\Bigr )\bigl (%
\sum_{t=1}^{3}D_{t}(e_{a})\bigr )\ -\ fh^{-1}(n+2)\,dh(e_{a})\bigl (%
\sum_{t=1}^{3}A_{t}(e_{a})\bigr ) \\
+\ {\nabla dh}(I_{1}e_{a},\xi _{1})\left (D_{1}(e_{a})\ -\ D_{2}(e_{a})\ -\
D_{3}(e_{a})\right )\ +\ {\nabla dh}(I_{2}e_{a},\xi _{2})\left
(-D_{1}(e_{a})\ +D_{2}(e_{a})\ -\ D_{3}(e_{a})\right) \\
+\ {\nabla dh}(I_{3}e_{a},\xi _{3})\left (-D_{1}(e_{a})\ -\ D_{2}(e_{a})\
+D_{3}(e_{a})\right ) \\
+\ h^{-1}\bigl (\,\sum_{s=1}^{3}dh(\xi _{s})\,dh(I_{s}e_{a})\bigr )\,\bigl (%
\sum_{t=1}^{3}D_{t}(e_{a})\bigr )\ +\ h^{-1}(n+2)\bigl (\sum_{s=1}^{3}\,dh(%
\xi _{s})\,dh(I_{s}e_{a})\bigr )\,\bigl (\sum_{t=1}^{3}A_{t}(e_{a})\bigr ) \\
+\ 4\sum_{s=1}^{3}{\nabla dh}(I_{s}e_{a},\xi _{s})\,A_{s}(e_{a})\ -\frac{10}{%
3}\ \bigl (\sum_{s=1}^{3}\ {\nabla dh}(I_{s}e_{a},\xi _{s})\bigr )\,\bigl (%
\sum_{t=1}^{3}A_{t}(e_{a})\bigr ),
\end{multline}%
where the last equality uses \eqref{e:F_s by D_s} to express the vectors $%
F_{s}$ by $D_{s}$, and the expansions of the vectors $A$ and $D$ according
to \eqref{d:A_s} and \eqref{d:def of D}. Since the dimension of $M$ is seven
it follows $U\ =\ \bar{U}\ =\ [\ {\nabla dh}-2h^{-1}dh\otimes dh]_{[3][0]}\
=\ 0$. This, together with the Yamabe equation
\eqref{e:conf
change scalar curv} , which when }$n=1$ becomes $\triangle
h=2-4h+3h^{-1}|\nabla h|^{2},${\ yield the formula, cf. \eqref{e:U by h},
\begin{equation}\label{e:1}
\nabla dh\,(X,\nabla h)\ +\ \sum_{s=1}^{3}\nabla dh\,(I_{s}X,I_{s}\nabla h)\
-\ \bigl (\ 2\ -\ 4h\ +3h^{-1}|\nabla h|^{2}\bigr )\,dh(X)\ =\ 0.
\end{equation}%
\noindent From equations \eqref{d:D_s} and \eqref{e:T--+} we have
\begin{equation*}
\begin{aligned} D_1(X)\ =\ h^{-2}\, dh(\xi_1)\,dh(I_1X)\ + \ \frac14h^{-2}\,
\bigl [\ \nabla dh\, (X,\nabla h)+\ \nabla dh\, (I_1X,I_1\nabla h)\\ - \
\nabla dh\, (I_2X,I_2\nabla h)\ -\ \nabla dh\, (I_3X,I_3\nabla h)\bigr ],\\
D_2(X)\ =\ h^{-2}\, dh(\xi_2)\, dh(I_2X)\  + \ \frac14h^{-2}\bigl[ \nabla
dh\, (X,\nabla h)\ -\ \nabla dh\, (I_1X,I_1\nabla h)\\ + \ \nabla dh\,
(I_2X,I_2\nabla h)\ -\ \nabla dh\, (I_3X,I_3\nabla h)\bigr],\\ D_3(X)\ =\
h^{-2}\, dh(\xi_3)\, dh(I_3X)\ +\ \frac14 h^{-2}\bigl[ \nabla dh(X,\nabla
h)-\ \nabla dh\, (I_1X,I_1\nabla h)\\ - \ \nabla dh\, (I_2X,I_2\nabla h)\ +
\ \nabla dh\, (I_3X,I_3\nabla h) \bigr]. \end{aligned}
\end{equation*}%
Expressing the first term in \eqref{e:1} by the rest and substituting with
the result in the above equations we come to%
\begin{multline}
D_{i}(e_{a})\ =\frac{1}{4}h^{-2}\left( 2\ -\ 4h\ +3h^{-1}|\nabla
h|^{2}\right) dh(e_{a})+h^{-2}\,dh(\xi _{i})\,dh(I_{i}e_{a})+ \\
\frac{1}{2}h^{-2}\left[ -\ \nabla dh\,(I_{j}e_{a},I_{j}\nabla h)\ -\nabla
dh\,(I_{k}e_{a},I_{k}\nabla h)\right] .  \label{e:D_s}
\end{multline}%
At this point, by a purely algebraic calculation, using
Lemma~\ref{l:transf of A} and \eqref{e:D_s} we find:
\begin{multline*}
\frac{22}{3}\,A_{1}\ -\ \frac{2}{3}\,A_{2}\ -\ \frac{2}{3}\,A_{3}\ +\
\frac{11}{3}\,\,D_{1}\ -\ \frac{1}{3}\,\,D_{2}\ -\ \frac{1}{3}\,\,D_{3} \\
=-3h^{-1}\left( 1+\frac{1}{2}h^{-1}dh(e_{a}) +\frac{1}{4}h^{-2}\lvert
\nabla h\rvert ^{2}\right) dh(e_{a})+3h^{-2}\left( \sum_{s=1}^{3}\,dh(\xi _{s})\,dh(I_{s}e_{a})\right)\\
+\frac{2}{3}h^{-1}{\nabla dh}(I_{1}e_{a},\xi _{1})-\frac{10}{3}h^{-1}{%
\nabla dh}(I_{2}e_{a},\xi _{2})-\frac{10}{3}h^{-1}{\nabla
dh}(I_{3}e_{a},\xi _{3}).
\end{multline*}
Similarly, %
\begin{multline*}
3\,A_{1}-A_{2}-A_{3}+2D_{1} =\left( -2h^{-1} + \frac{1}{2}h^{-2}
+h^{-3}|\nabla h|^{2}\right) dh(e_{a})-\frac{1}{2}h^{-2}\sum_{s=1}^3{\nabla dh}(I_{s}e_{a},I_{s}\nabla h) \\
+h^{-1}\ {\nabla dh}(I_{1}e_{a},\xi _{1})-h^{-1}\ {\nabla dh}%
(I_{2}e_{a},\xi _{2})-h^{-1}\ {\nabla dh}(I_{3}e_{a},\xi _{3})
+h^{-2}\sum_{s=1}^3dh(\xi _{s})\,dh(I_{s}e_{a}).
\end{multline*}
On the other hand, the coefficient of {$A_{1}(e_{a})$ in \eqref{e:diver}
is found to be, after setting }$\ n=1$,
\begin{multline*}
h\Big[ -\ 3\left( \ 1\ +\frac{1}{2}h^{-1}\ +\ \frac{1}{4}h^{-2}\lvert
\nabla h\rvert ^{2}\right) h^{-1}\,dh(e_{a})+\ 3h^{-2}\left(
\sum_{s=1}^{3}\,dh(\xi _{s})\,dh(I_{s}e_{a})\right)\\
 +\
\frac{2}{3}h^{-1}{\nabla dh}(I_{1}e_{a},\xi _{1})-\frac{10}{3}h^{-1}{
\nabla dh}(I_{2}e_{a},\xi _{2})-\frac{10}{3}h^{-1}{\nabla
dh}(I_{3}e_{a},\xi _{3}) \Big],
\end{multline*}
while the coefficient of $D_{1}(e_{a})$ {\ in \eqref{e:diver} is}
\begin{multline}
dh(e_{a})\ -\ \frac{1}{4}h^{-2}dh(e_{a})\lvert \nabla h\rvert ^{2}\ +\ \frac{%
1}{2}h^{-1}\,\ {\nabla dh}(e_{a},\nabla h)-\ f\,h^{-1}dh(e_{a})\  \\
+\ {\nabla dh}(I_{1}e_{a},\xi _{1})\ -\ {\nabla dh}(I_{2}e_{a},\xi _{2})-\ {%
\nabla dh}(I_{3}e_{a},\xi _{3})D_{1}(e_{a}) \\
+\ h^{-1}\bigl (\,\sum_{s=1}^{3}dh(\xi _{s})\,dh(I_{s}e_{a})\bigr ).
\end{multline}%
Substituting ${\nabla dh}(e_{a},\nabla h)$ according to \eqref{e:1}, i.e.,
$\nabla dh\,(e_{a},\nabla h)\ =-\ \sum_{s=1}^{3}\nabla
dh\,(I_{s}e_{a},I_{s}\nabla h)\ +\ \bigl (\ 2\ -\ 4h\ +3h^{-1}|\nabla
h|^{2}\bigr )\,dh(e_{a})$ and using the definition of $f$ transforms the
above expression into
\begin{multline*}
dh(e_{a})\ -\ \frac{1}{4}h^{-2}dh(e_{a})\lvert \nabla h\rvert ^{2}\ -\
\left( \frac{1}{2}\ +\ h\ +\ \frac{1}{4}h^{-1}\lvert \nabla h\rvert
^{2}\right) \,h^{-1}dh(e_{a})\  \\
+\ \frac{1}{2}h^{-1}\,\ \left( -\ \sum_{s=1}^{3}\nabla
dh\,(I_{s}e_{a},I_{s}\nabla h)\ +\ \bigl (\ 2\ -\ 4h\ +3h^{-1}|\nabla h|^{2}%
\bigr )\,dh(e_{a})\right) \\
+\ {\nabla dh}(I_{1}e_{a},\xi _{1})\ -\ {\nabla dh}(I_{2}e_{a},\xi _{2})-\
{\nabla dh}(I_{3}e_{a},\xi _{3})D_{1}(e_{a})+\ h^{-1}\,\left(
\sum_{s=1}^{3}dh(\xi _{s})\,dh(I_{s}e_{a})\right).
\end{multline*}
Simplifying the above expression shows that the coefficient of $D_{1}(e_{a})$ {\ in \eqref{e:diver} is}%
\begin{multline*}
\left( -2\ +\frac{1}{2}h^{-1}+h^{-2}|\nabla h|^{2}\right) dh(e_{a})-\
\frac{1}{2}h^{-1}\left( \sum_{s=1}^{3}\nabla dh\,(I_{s}e_{a},I_{s}\nabla
h)\
\right)  \\
+\ {\nabla dh}(I_{1}e_{a},\xi _{1})\ -\ {\nabla dh}(I_{2}e_{a},\xi _{2})-\
{\nabla dh}(I_{3}e_{a},\xi _{3})+\ h^{-1}\,\left( \sum_{s=1}^{3}dh(\xi
_{s})\,dh(I_{s}e_{a})\right)
\end{multline*}

Hence, we proved that the coefficient of $D_{1}(e_{a})$  in
\eqref{e:diver} is $h\left(
3\,A_{1}-A_{2}-A_{3}+2D_{1}\right)(e_a) $, while
those of $A_{1}(e_{a})$ is  $h\bigl (\frac{22}{3}%
\,A_{1}\ -\ \frac{2}{3}\,A_{2}\ -\ \frac{2}{3}\,A_{3}\ +\ \frac{11}{3}%
\,\,D_{1}\ -\ \frac{1}{3}\,\,D_{2}\ -\ \frac{1}{3}\,\,D_{3}\bigr
)(e_{a})$.  A cyclic permutation gives the rest of the coefficients in
\eqref{e:diver}. With this, the divergence \eqref{e:diver} can be written
in the form
\begin{multline*}
\nabla ^{\ast }\Bigl(fD\ +\ \sum_{s=1}^{3}dh(\xi _{s})\,F_{s}\ +\
4\sum_{s=1}^{3}dh(\xi _{s})I_{s}A_{s}\ -\ \frac{10}{3}\sum_{s=1}^{3}dh(\xi
_{s})\,I_{s}A\Bigr) \\
=\ f\lvert T^{0}\rvert ^{2}\ +\ h\,\sigma _{1,2,3}\Bigl \{\ g\,\bigl(%
\,D_{1},\,3A_{1}\ -\ A_{2}\ -\ A_{3}\ +\ 2\,D_{1}\bigr) \\
+\ g\,\bigl(\,A_{1},\,\frac{22}{3}\,A_{1}\ -\ \frac{2}{3}\,A_{2}\ -\ \frac{2%
}{3}\,A_{3}\ +\ \frac{11}{3}\,\,D_{1}\ -\ \frac{1}{3}\,\,D_{2}\ -\ \frac{1}{3%
}\,\,D_{3}\bigr)\Bigr \},
\end{multline*}%
\noindent where $\sigma _{1,2,3}$ denotes the sum over all positive
permutations of $(1,2,3)$.} 
Let $Q$ be equal to
\begin{equation*}
Q:=\left[ {%
\begin{array}{cccccc}
2\, & 0 & 0 & {\displaystyle\frac{10}{3}}\, & -{\displaystyle\frac{2}{3}}\,
& -{\displaystyle\frac{2}{3}}\, \\[2ex]
0 & 2\, & 0 & -{\displaystyle\frac{2}{3}}\, & {\displaystyle\frac{10}{3}}\,
& -{\displaystyle\frac{2}{3}}\, \\[2ex]
0 & 0 & 2\, & -{\displaystyle\frac{2}{3}}\, & -{\displaystyle\frac{2}{3}}\,
& {\displaystyle\frac{10}{3}}\, \\[2ex]
{\displaystyle\frac{10}{3}}\, & -{\displaystyle\frac{2}{3}}\, & -{%
\displaystyle\frac{2}{3}}\, & {\displaystyle\frac{22}{3}}\, & -{\displaystyle%
\frac{2}{3}}\, & -{\displaystyle\frac{2}{3}}\, \\[2ex]
-{\displaystyle\frac{2}{3}}\, & {\displaystyle\frac{10}{3}}\, & -{%
\displaystyle\frac{2}{3}}\, & -{\displaystyle\frac{2}{3}}\, & {\displaystyle%
\frac{22}{3}}\, & -{\displaystyle\frac{2}{3}}\, \\[2ex]
-{\displaystyle\frac{2}{3}}\, & -{\displaystyle\frac{2}{3}}\, & {%
\displaystyle\frac{10}{3}}\, & -{\displaystyle\frac{2}{3}}\, & -{%
\displaystyle\frac{2}{3}}\, & {\displaystyle\frac{22}{3}}\,%
\end{array}%
}\right]
\end{equation*}%
\noindent so that
\begin{equation*}
\nabla ^{\ast }\Bigl(fD\ +\ \sum_{s=1}^{3}dh(\xi _{s})\,F_{s}\ +\ 4dh(\xi
_{s})\,I_{s}A_{s}\ -\ \frac{10}{3}\sum_{s=1}^{3}dh(\xi _{s})\,I_{s}A\Bigr)\
=\ f\,\lvert T^{0}\rvert ^{2}\ +\ h\,\langle QV,\,V\rangle ,
\end{equation*}%
with $V=(D_{1},D_{2},D_{3},A_{1},A_{2},A_{3})$. It is not hard to see that
the eigenvalues of $Q$ are given by
\begin{equation*}
\{0,\quad 0,\quad 2\,(2+\sqrt{2}),\quad 2\,(2-\sqrt{2}),\quad 10,\quad 10\},
\end{equation*}%
\noindent which shows that $Q$ is a non-negative matrix.
\end{proof}


\section{Proofs of the main theorems}

The proofs rely on Theorem~\ref{t:div formula} and the following
characterization of all qc-Einstein structures conformal to the standard qc
structures on the Heisenberg group.

\begin{thrm}
\cite[Theorem~1.2]{IMV}\label{nnn} Let $\Theta=\frac{1}{2h}\tilde\Theta$ be
a conformal deformation of the standard qc-structure $\tilde\Theta$ on the
quaternionic Heisenberg group $\boldsymbol{G\,(\mathbb{H})}$. Then $\Theta$
is qc-Einstein if and only if, up to a left translation the function $h$ is
given by
\begin{equation}  \label{e:Yamabe solutions}
h \ =\ c\ \Big [ \big ( 1\ +\ \nu\, |q|^2 \big )^2\ +\ \nu^2\, (x^2\ +\ y^2\
+\ z^2)\Big ],
\end{equation}
where $c$ and $\nu$ are any positive constants.
\end{thrm}

Consider first the case of the (seven dimensional) sphere.


\subsection{Proof of Theorem 1.1}

Integrating the divergence formula of Theorem~\ref{t:div formula} we see
that according to the divergence theorem established in \cite[Proposition 8.1%
]{IMV} the integral of the left-hand side is zero. Thus, the right-hand side
vanishes as well, which shows that the quaternionic contact structure $\eta$
has vanishing torsion, i.e., it is also qc-Einstein {\ according to Theorem~%
\ref{sixtyseven}.}

Next we bring into consideration the 7-dimensional quaternionic Heisenberg
group and the quaternionic Cayley transform as described in \cite[Section 5.2%
]{IMV}. The quaternionic Heisenberg group of dimension $7$ is $\boldsymbol{%
G\,(\mathbb{H})} \ =\mathbb{H}\times\text {Im}\, \mathbb{H}$. The group law
is given by \hspace{2mm} $( q^{\prime }, \omega^{\prime })\ =\ (q_o,
\omega_o)\circ(q, \omega)\ =\ (q_o\ +\ q, \omega\ +\ \omega_o\ + \ 2\ \text {%
Im}\ q_o\, \bar q), $ \noindent where $q,\ q_o\in\mathbb{H}$ and $\omega,
\omega_o\in \text {Im}\, \mathbb{H}$. The left-invariant orthonormal basis
of the horizontal space is
\begin{equation*}
\begin{aligned} T_{1}\  =\  {\frac{\partial }{\partial t_1}} {}\
+2x^{1}\frac{\partial {}}{\partial x}+2y^{1}\frac{\partial {}}{\partial
y}+2z^{1}\frac{\partial {}}{\partial z} \,,\quad X_{1}\  =\  {\frac{\partial
}{\partial x_1}} {}\ -2t^{1}\frac{\partial {}}{\partial
x}-2z^{1}\frac{\partial {}}{\partial y}+2y^{1}\frac{\partial {}}{\partial z}
\,\\ Y_{1}\  =\  {\frac{\partial }{\partial y_1}} {}\ +2z^{1}\frac{\partial
{}}{\partial x}-2t^{1}\frac{\partial {}}{\partial y}-2x^{1}\frac{\partial
{}}{\partial z}\,,\quad Z_{1}\  =\  {\frac{\partial }{\partial z_1}} {}\
-2y^{1}\frac{\partial {}}{\partial x}+2x^{1}\frac{\partial {}}{\partial
y}-2t^{1}\frac{\partial {}}{\partial z}\, \end{aligned}
\end{equation*}
\noindent using $q\ =\ t_1\ +\ i\, x_1 \ +\ j\, y_1\ +\ k\, z_1$ and $\omega
\ =\ i\, x \ +\ j\, y\ +\ k\, z$. The central (vertical) orthonormal vector
fields $\xi_1,\xi_2,\xi_3$ are described as follows
\begin{equation*}
\xi_1=2\frac{\partial {}}{\partial x}\,\qquad \xi_2=2\frac{\partial {}}{%
\partial y}\,\qquad \xi_3=2\frac{\partial {}}{\partial z}\,.
\end{equation*}
\noindent {\ Let us identify the (seven dimensional) group $\boldsymbol{G\,(%
\mathbb{H})}$ with the boundary $\Sigma$ of a Siegel domain in $\mathbb{H}%
\times\mathbb{H}$,}
\begin{equation*}
\Sigma\ =\ \{ (q^{\prime },p^{\prime })\in \mathbb{H}\times\mathbb{H}\ :\
\Re {\ p^{\prime }}\ =\ \lvert q^{\prime 2 }\}.
\end{equation*}
$\Sigma$ carries a natural group structure and the map $(q, \omega)\mapsto
(q,\lvert q \rvert^2 - \omega)\,\in\,\Sigma$ is an isomorphism between $%
\boldsymbol{G\,(\mathbb{H})}$ and $\Sigma$.

{\ The standard contact form, written as a purely imaginary quaternion
valued form, on $\boldsymbol{G\,(\mathbb{H})}$ is given by \hspace{2mm} $2%
\tilde{\Theta} = (d\omega \ - \ q \cdot d\bar q \ + \ dq\, \cdot\bar q),$
\hspace{2mm} where $\cdot$ denotes the quaternion multiplication. Since
\hspace{2mm} $dp\ =\ q\cdot d\bar q\ +\ dq\, \cdot\bar {q}\ -\ d\omega, $
\hspace{2mm} under the identification of $\boldsymbol{G\,(\mathbb{H})}$ with
$\Sigma$ we also have \hspace{2mm} $2\tilde{\Theta}\ =\ - dp^{\prime }\ +\
2dq^{\prime }\cdot\bar {q}^{\prime }. $ \hspace{2mm} Taking into account
that $\tilde{\Theta}$ is purely imaginary, the last equation can be written
also in the following form
\begin{equation*}
4\,\tilde{\Theta}\ =\ (d\bar p^{\prime }\ -\ d p^{\prime })\ +\ 2dq^{\prime
}\cdot\bar {q^{\prime }}\ -\ 2q^{\prime }\cdot d\bar q^{\prime }.
\end{equation*}
} The (quaternionic) Cayley transform is the map $\mathcal{C}:\,
S\setminus\{(-1,0)\}\,\mapsto \Sigma$ from the sphere $S\ =\ \{(q,p)\in
\mathbb{H}\times\mathbb{H}\, :\ \lvert q \rvert^2+\lvert p \rvert^2=1
\}\subset \mathbb{H}\times\mathbb{H}$ minus a point to the Heisenberg group $%
\Sigma= \{ (q_1,p_1)\in \mathbb{H}\times\mathbb{H}\ :\ \Re {\ p_1}\ =\
\lvert q_1 \rvert^2 \}$, with $\mathcal{C}$ defined by
\begin{equation}  \label{d:Cayley}
(q_1, p_1)\ =\ \mathcal{C}\ \Big ((q, p)\Big),\qquad q_1\ =\ (1+p)^{-1} \
q,\quad p_1\ =\ (1+p)^{-1} \ (1-p).
\end{equation}

\noindent with an inverse $(q, p)\ =\ \mathcal{C}^{-1}\Big ((q_1, p_1)\Big)$
given by
\begin{equation}  \label{d:inverse Cayley}
q\ =\ \ 2 \,(1+p_1)^{-1}\,q_1, \quad p\ =\ (1-p_1)\,(1+p_1)^{-1}.
\end{equation}

\noindent The Cayley transform is a conformal quaternionic contact
diffeomorphism between the quaternionic Heisenberg group with its standard
quaternionic contact structure $\tilde\Theta$ and $S\setminus\{(-1,0)\}$
with its standard structure $\tilde\eta$, see \cite{IMV},
\begin{equation}  \label{e:Cayley transf of ctct form}
\lambda\ \cdot (\mathcal{C}_*\, \tilde\eta)\ \cdot \bar\lambda\ =\ \frac {8}{%
\lvert 1+p_1\, \rvert^2}\, \tilde\Theta,
\end{equation}
where $\lambda\ =\ \frac {1+p_1}{\lvert 1+p_1\, \rvert}$ is a unit
quaternion and $\tilde\eta$ is the standard quaternionic contact form on the
sphere, \hspace{5mm} $\tilde\eta\ =\ dq\cdot \bar q\ +\ dp\cdot \bar p\ -\
q\cdot d\bar q -\ p\cdot d\bar p. $ Hence, up to a constant multiplicative
factor and a quaternionic contact automorphism the forms $\mathcal{C}%
_*\tilde\eta$ and $\tilde\Theta$ are conformal to each other. It follows
that the same is true for $\mathcal{C}_*\eta$ and $\tilde\Theta$. In
addition, $\tilde\Theta$ is qc-Einstein by definition, while $\eta$ and
hence also $\mathcal{C}_* \eta$ are qc-Einstein as we observed at the
beginning of the proof. According to Theorem~\ref{nnn}, up to a
multiplicative constant factor, the forms $\mathcal{C}_*\tilde\eta$ and $%
\mathcal{C}_*\eta$ are related by a translation or dilation on the
Heisenberg group. Hence, we conclude that up to a multiplicative constant, $%
\eta$ is obtained from $\tilde\eta$ by a conformal quaternionic contact
automorpism which proves the first claim of Theorem \ref{t:Yamabe}. From the
conformal properties of the Cayley transform and \cite{Va1,Va2} it follows
that the minimum $\lambda(S^{4n+3})$ is achieved by a smooth 3-contact form,
{\ which due to the Yamabe equation is of constant qc-scalar curvature.}
This shows the second claim of Theorem \ref{t:Yamabe}. 

\medskip 

\subsection{Proof of Theorem 1.3}

{Let $\mathcal{D}^{1,2}$ be the space of functions $u\in L^{2^*}(\boldsymbol{%
G\,(\mathbb{H})})$ having distributional horizontal gradient $|\nabla
u|^2=|T_1u|^2+|X_1u|^2+|Y_1u|^2+|Z_1u|^2\in L^2(\boldsymbol{G\,(\mathbb{H})}%
) $ with respect to the Lebesgue measure $dH$ on $\mathbb{R}^7$, which is
the Haar measure on the group. Let us define the constant ($2^*=5/2$ here)
\begin{equation*}  \label{S}
\Lambda\ \overset{def}{=}\ inf\ \left\{\underset{\boldsymbol{G\,(\mathbb{H})}%
}{\int} |\nabla v|^2\ dH\ : \ v\in \mathcal{D}^{1,2},\ v \geq 0,\ \underset{%
\boldsymbol{G\,(\mathbb{H})}}{\int} |v|^{2^*}\ dH\ =\ 1 \right\}.
\end{equation*}
Let $v$ be a function for which the infimum is achieved. Note that such
function exists by \cite{Va1} or \cite{Va2}. Furthermore, $\Lambda\ {=}\
S_2^{-2}$, where $S_2$ is the best constant in the $L^2$ Folland-Stein
inequality \eqref{FS}, since $v\in \mathcal{D}^{1,2}$ implies $|v|\in
\mathcal{D}^{1,2}$ and the gradient is the same a.e.. From the choice of $v$
we have
\begin{equation*}
\Lambda\ =\ \underset{\boldsymbol{G\,(\mathbb{H})}}{\int}\ |\nabla v|^2\
dH,\quad\quad\quad \underset{\boldsymbol{G\,(\mathbb{H})}}{\int}\ v^{2^*}\
dH\ =\ 1.
\end{equation*}
Writing the Euler-Lagrange equation of the constrained problem we see that $%
v $ is a non-negative entire solution of $\big ( T^2_1 \ +\ X^2_1 \ +\ Y^2_1
\ +\ Z^2_1 \big ) v = - \Lambda\ v^{3/2}$. By \cite[Lemma 10.2]{GV2} ( see
\cite{Va2} or \cite[Theorem 10.3]{Va1} for further details) $v$ is a bounded
function. Similarly to \cite[Theorem 16.7]{FS} it follows $v$ is a Lipschits
continuous function in the sense of non-isotropic Lipschits spaces \cite{F}.
Iterating this argument and using \cite[Theorem 5.25]{F} we see that $v$ is
a $C^\infty$ smooth function on the set where it is positive, while being of
class $\Gamma_{loc}^{2,\beta}$, the non-isotropic Lipschits space, for some $%
\beta>0$. In particular $v$ is continuously differentiable function by \cite[%
Theorem 5.25]{F}. Applying the Hopf lemma \cite[Theorem 2.13]{GV} on the set
where $v$ is positive shows that $v$ cannot vanish, i.e., it is a positive
entire solution to the Yamabe equation. The positivity can also be seen by
the Harnack inequality, see \cite{Wei} for example. } Let $u\ \overset{def}{=%
}\ \Lambda^{\frac {1}{2\text{*}-2}}\ v$, then $u$ is a positive entire
solution of the Yamabe equation
\begin{equation}  \label{e:Yamabe}
\big ( T^2_1 \ +\ X^2_1 \ +\ Y^2_1 \ +\ Z^2_1 \big ) \, u\ = \ - u^{3/2}
\end{equation}
>From the definition of $u$, we have
\begin{equation*}
\Lambda\ = \ \Bigl(\underset{\boldsymbol{G\,(\mathbb{H})}}{\int} |\nabla
u|^2\ dH\Bigr)^{\frac {1}{5}}\ =\ \Bigl( \underset{\boldsymbol{G\,(\mathbb{H}%
)}}{\int}\ u^{5/2}\ dH\Bigr)^{\frac {1}{5}}.
\end{equation*}
We shall compute the last integral by determining $u$ with the help of the
divergence formula.

{As before, let $\tilde\Theta$ be the standard contact form on $\boldsymbol{%
G\,(\mathbb{H})}$ identified with $\Sigma$. Using the inversion and the
Kelvin transform on $\boldsymbol{G\,(\mathbb{H})}$, cf. \cite[Sections 8 and
9]{GV2}, we can see that if $\Theta=\frac {1}{2h}\tilde\Theta$ has constant
scalar curvature, then the Cayley transform lifts the qc structure defined
by $\Theta$ to a qc structure of constant qc-scalar curvature on the sphere,
which is conformal to the standard.} The details are as follows. Let us
define two contact forms $\Theta_1$ and $\Theta_2$ on $\Sigma$ setting
\begin{equation*}
\Theta_1\ =\ u^{4/(Q-2)}\tilde\Theta, \quad \text{and}\quad \Theta_2\ =\ ({%
\mathcal{K}u})^{4/(Q-2)}\ \frac {\bar p^{\prime }}{|p^{\prime }|}\
\tilde\Theta\ \frac {p^{\prime }}{|p^{\prime }|},
\end{equation*}
where $u$ is as in \eqref{e:Yamabe}, $\mathcal{K}u$ is its Kelvin transform,
see \eqref{e:Kelvin} for the exact formula, and $Q$ is the homogeneous
dimension of the group. Notice that $\frac {\bar p^{\prime }}{|p^{\prime }|}%
\ \tilde\Theta\ \frac {p^{\prime }}{|p^{\prime }|}$ defines the same qc
structure on the Heisenberg group as $\tilde\Theta$ and $\mathcal{K}u$ is a
smooth function on the whole group according to \cite[Theorem 9.2]{GV2}. We
are going to see that using the Cayley transform these two contact forms
define a contact form on the sphere, which is conformal to the standard and
has constant qc-scalar curvature.

Let $P_1=(-1,0)$ and $P_2=(1,0)$ be correspondingly the 'south' and 'north'
poles of the unit sphere $S\ =\ \{\lvert q \rvert^2+\lvert p \rvert^2=1 \}$.
Let $\mathcal{C}_1$ and $\mathcal{C}_2$ be the corresponding Cayley
transforms defined, respectively, on $S\setminus\{P_1\}$ and $%
S\setminus\{P_2\}$. Note that $\mathcal{C}_1$ was defined in \eqref{d:Cayley}%
, while $\mathcal{C}_2$ is given by

\begin{equation}  \label{d:2nd Cayley}
(q_2, p_2)\ =\ \mathcal{C}_2\ \Big ((q, p)\Big),\qquad q_2\ =\ -(1-p)^{-1} \
q,\quad p_2\ =\ (1-p)^{-1} \ (1+p).
\end{equation}
In order that $\Theta_1$ and $\Theta_2$ define a contact form $\eta$ on the
sphere it is enough to see that
\begin{equation}  \label{e:global form}
\Theta_1 (p,q)\ =\ \Theta_2 \circ\mathcal{C}_2\circ\mathcal{C}_1^{-1} (p,q),
\quad \text{i.e.,}\quad \Theta_1 \ =\ (\mathcal{C}_2\circ\mathcal{C}%
_1^{-1})^* \Theta_2.
\end{equation}
A calculation shows that $\mathcal{C}_2\circ\mathcal{C}_1^{-1}:\Sigma%
\rightarrow\Sigma$ is given by
\begin{equation*}
q_2\ =\ -p_1^{-1}\, q_1,\qquad p_2\ =\ p_1^{-1},
\end{equation*}
or, equivalently, in the model $\boldsymbol{G\,(\mathbb{H})}$
\begin{equation*}
q_2\ =\ -(|q_1|^2-\omega_1)^{-1}\, q_1, \qquad \omega_2\ =\ -\frac {\omega_1%
}{|q_1|^4+|\omega_1|^2}.
\end{equation*}
Hence, $\sigma\ =\ \mathcal{C}_2\circ\mathcal{C}_1^{-1} $ is {an involution}
on the group. Furthermore, with the help of
\eqref{e:Cayley transf of ctct
form} we calculate
\begin{equation*}
{\mathcal{C}_1}_*\circ{\mathcal{C}_2}^*\ \Theta\ =\ \frac {1}{|p_1|^2}
\,\bar\mu\, \Theta\, \mu, \qquad \mu\ =\ \frac {p_1}{|p_1|},
\end{equation*}

\noindent which proves the identity \eqref{e:global form}. Using the
properties of the Kelvin transform, \cite[Sections 8 and 9]{GV2},
\begin{equation}  \label{e:Kelvin}
(\mathcal{K}u)\, (q^{\prime },p^{\prime })\ \overset{def}{=}\ |p^{\prime
-(Q-2)/2}\, u\bigl (\sigma(q^{\prime },p^{\prime })\bigr ),
\end{equation}
we see that $u$ and $\mathcal{K}u$ are solutions of the Yamabe equation %
\eqref{e:Yamabe}. This implies that the contact form $\eta$ has constant
qc-scalar curvature, equal to $\frac {4(Q+2)}{Q-2}.$

Notice that $\eta$ is conformal to the standard form $\tilde\eta$ and the
arguments in the preceding proof imply then that $\eta$ is qc-Einstein. A
small calculation shows that this is equivalent to the fact that if we set

\begin{equation}  \label{e:Yamabe sols}
\bar u\ =\ 2^{10}\ [(1+\lvert q \rvert^2)^2\ +\ \lvert \omega \rvert^2]^{-2},
\end{equation}
\noindent then $\bar u$ satisfies the Yamabe equation \eqref{e:Yamabe} and
all other nonnegative solutions of \eqref{e:Yamabe} in the space $\mathcal{D}%
^{1,2}$ are obtained from $\bar u$ by translations and dilations,
\begin{gather} \label{e:translation}
\tau_{(q_o,\omega_o)} \bar u\ (q,\omega)\ \overset{def}{=}\ \bar u (q_o+q,
\omega + \omega_o), \\ \label{e:scaling}\bar u_\lambda\, (q)\
\overset{def}{=}\ \lambda^{4}\bar u(\lambda q,\lambda^2 \omega),
\quad\quad\quad \lambda >0.
\end{gather}

Thus, $u$ which was defined in the beginning of the proof is given by
equation \eqref{e:Yamabe sols} up to translations and dilations. This allows
the calculation of the best constant in the Folland-Stein inequality, see
\cite[(4.52)]{GV},
\begin{equation*}
\Lambda ^5\ =\ \underset{\boldsymbol{G\,(\mathbb{H})}}{\int}\ \frac{2^{25}}{%
\left[(1 + |q|^2)^2\ +\ |\omega|^2)\right]^{5}}\ dH\ =\ 2^{25}\,\pi^{7/2}
\frac{\Gamma(\frac{7}{2})}{\Gamma(7)}\ =\ \frac {\pi^{12/10}}{12},
\end{equation*}
where $\Gamma$ is the Gamma function. Hence
\begin{equation*}
S_2\ =\ \Lambda^{-1/2}\ =\ \frac{2\sqrt{3}}{\pi^{3/5}}.
\end{equation*}
Recalling the relation between $u$ and $v$ we find that the extremals in the
Folland-Stein embedding are given by
\begin{equation*}
v \ =\ \frac{2^{11}\sqrt{3}}{\pi^{3/5}}[(1+\lvert q \rvert^2)^2\ +\ \lvert
\omega \rvert^2]^{-2}
\end{equation*}
and its translations and dilations. The proof of Theorem \ref{t:FS} is
complete.



\begin{thebibliography}{Biq1}
\bibitem[AK]{AK} Alekseevsky, D. \& Kamishima, Y., \emph{Pseudo-conformal
quaternionic CR structure on $(4n+3)$-dimensional manifold}, math.GT/0502531.

\bibitem[Biq1]{Biq1} Biquard, O., \emph{M\'{e}triques d'Einstein
asymptotiquement sym\'{e}triques}, Ast\'{e}risque \textbf{265} (2000).

\bibitem[Biq2]{Biq2} Biquard, O., \emph{Quaternionic contact structures},
Quaternionic structures in mathematics and physics (Rome, 1999), 23--30
(electronic), Univ. Studi Roma "La Sapienza", Roma, 1999.

\bibitem[CSal]{CSal} Capria, M. \& Salamon, S., \emph{Yang-Mills fields on
quaternionic spaces} Nonlinearity \textbf{1} (1988), no. 4, 517--530.

\bibitem[D]{D} Duchemin, D., \emph{Quaternionic contact structures in
dimension 7}, Ann. Inst. Fourier (Grenoble) \textbf{56} (2006), no. 4,
851--885.

\bibitem[D1]{D1} \bysame, \emph{Quaternionic contact hypersurfaces},
math.DG/0604147.

\bibitem[F]{F} Folland, G., \emph{Subelliptic estimates and function spaces
on nilpotent Lie groups}, Ark. Math., \textbf{13} (1975), 161--207.


\bibitem[FSt]{FS} Folland, G. B. \& Stein, E. M., \emph{Estimates for the $%
\Bar {\partial}_{b}$ Complex and Analysis on the Heisenberg Group}, Comm.
Pure Appl. Math., \textbf{27}~(1974), 429--522.


\bibitem[GV1]{GV} Garofalo, N. \& Vassilev, D.,\ \emph{\ Symmetry properties
of positive entire solutions of Yamabe type equations on groups of
Heisenberg type}, Duke Math J, \textbf{106} (2001), no. 3, 411--449.

\bibitem[GV2]{GV2} {\leavevmode\vrule height 2pt depth -1.6pt width 23pt},\
\textit{Regularity near the characteristic set in the non-linear Dirichlet
problem and conformal geometry of sub-Laplacians on Carnot groups}, Math
Ann., \textbf{318} (2000), no. 3, 453--516.

\bibitem[IMV]{IMV} Ivanov, St., Minchev, I., \& Vassilev, D.,\ \emph{%
Quaternionic contact Einstein structures and the quaternionic contact Yamabe
problem}, preprint, math.DG/0611658.

\bibitem[JL1]{JL2} Jerison, D., \& Lee, J., \emph{The Yamabe problem on $CR$
manifolds}, J. Diff. Geom., \textbf{25}~(1987), 167--197.

\bibitem[JL2]{JL3} \bysame, \emph{Extremals for the Sobolev inequality on
the Heisenberg group and the CR Yamabe problem}, J. Amer. Math. Soc.,
\textbf{1}~(1988), no. 1, 1--13.

\bibitem[LP]{LP} Lee, J. M. \& Parker, T., \emph{\ The Yamabe Problem}, Bull
Am. Math. Soc. \textbf{17}~(1987), no. 1, 37--91.

\bibitem[M]{M} Mostow, G. D., \emph{Strong rigidity of locally symmetric
spaces}, Annals of Mathematics Studies, No. 78. Princeton University Press,
Princeton, N.J.; University of Tokyo Press, Tokyo, 1973. v+195 pp.

\bibitem[Ob]{Ob} Obata, M., \emph{The conjecture of conformal
transformations in Riemannian manifolds}, J. Diff. Geom., \textbf{6} (1971),
247--258.

\bibitem[P]{P} Pansu, P., \emph{M\'etriques de Carnot-Carath\'eodory et
quasiisom\'etries des espaces sym\'etriques de rang un}, Ann. of Math. (2)
129 (1989), no. 1, 1--60.

\bibitem[Va]{Va2} Vassilev, D, \textit{Yamabe type equations on Carnot groups%
}, Ph. D. thesis Purdue University, 2000.

\bibitem[Va1]{Va1} Vassilev, D., \textit{Regularity near the characteristic
boundary for sub-laplacian operators}, Pacific J Math, \textbf{227} (2006),
no.~2, 361--397.

\bibitem[W]{Wei} Wang, W., \emph{The Yamabe problem on quaternionic contact
manifolds}, Ann. Mat. Pura Appl., \textbf{186} (2007), no. 2, 359--380.

\end{thebibliography}
\end{document}